\documentclass[11pt]{amsart}
\pdfoutput=1
\usepackage{graphicx}
\usepackage{epsfig}
\usepackage[latin1]{inputenc}
\usepackage{amsmath}
\usepackage{amsfonts}
\usepackage{amssymb}
\usepackage{amsthm}
\usepackage{amscd}
\usepackage{verbatim}
\usepackage{subfigure}
\usepackage{pinlabel}
\usepackage{stmaryrd}

\usepackage{tikz}
\usetikzlibrary{arrows}
\usepackage{verbatim}
	
\newtheorem{theorem}{Theorem}[section]
\newtheorem{Theorem}{Theorem}
\newtheorem{lemma}[theorem]{Lemma}
\newtheorem{proposition}[theorem]{Proposition}
\newtheorem{Proposition}[Theorem]{Proposition}

\newtheorem{definition}[theorem]{Definition}
\newtheorem{Definition}[Theorem]{Definition}
\newtheorem{remark}[theorem]{Remark}

\def\Z{\mathbb{Z}}
\def\N{\mathbb{N}}

\def\F{\mathbb{F}}

\def\cF{\mathcal{F}}
\def\cC{\mathcal{C}}

\def\CFKi{CFK^{\infty}}
\def\horz{\textup{horz}}
\def\vert{\textup{vert}}

\title{The knot Floer complex and the smooth concordance group}

\subjclass[2009]{}
\author{Jennifer Hom}

\address{Department of Mathematics, Columbia University, New York, NY 10027
\newline\indent{\tt hom@math.columbia.edu}}

\numberwithin{equation}{section}

\begin{document}

\begin{abstract}
We define a new smooth concordance homomorphism based on the knot Floer complex and an associated concordance invariant, $\varepsilon$. As an application, we show that an infinite family of topologically slice knots are independent in the smooth concordance group.
\end{abstract}

\maketitle

\section{Introduction}

The set of isotopy classes of knots in $S^3$, under the operation of connected sum, forms a monoid.
Two knots are \emph{concordant} if they cobound a smooth, properly embedded cylinder in $S^3 \times [0,1]$. The monoid of knots, modulo concordance, forms the \emph{concordance group}, denoted $\cC$. If we loosen the conditions and only require that the cylinder be locally flat, rather than smooth, we obtain the \emph{topological concordance group}. Understanding the difference between these two groups sheds some light on the distinction between the smooth and topological categories.

Ozsv\'ath and Szab\'o \cite{OSknots}, and independently Rasmussen \cite{R}, defined an invariant, knot Floer homology, associated to a knot in $S^3$. This invariant comes in many different flavors, the most robust being $CFK^{\infty}(K)$, a $\Z$-filtered chain complex over the ring $\F[U, U^{-1}]$, where $\F=\Z/2\Z$ and $U$ is a formal variable. There is a second filtration induced by $-$($U$-exponent) allowing us to view $\CFKi(K)$ as a $\Z \oplus \Z$-filtered chain complex. The filtered chain homotopy type of this complex is an invariant of the knot $K$. The weaker invariant, $\widehat{CFK}(K)$, takes the form of a $\Z$-filtered chain complex over $\F$, and is obtained by taking the degree zero part of the associated graded object with respect to one of the filtrations.

Within the complex $\widehat{CFK}(K)$ lives a $\Z$-valued concordance invariant, $\tau(K)$, defined by Ozsv\'ath and Szab\'o in \cite{OS4ball}. The total homology of $\widehat{CFK}(K)$ has rank one, and $\tau$ measures the minimum filtration level where this homology is supported. The invariant $\tau$ gives a surjective homomorphism from the smooth concordance group $\cC$ to the integers:
$$\tau: \mathcal{C} \rightarrow \Z,$$
which gives a new proof of the Milnor conjecture \cite{OS4ball} and is strong enough to obstruct topologically slice knots from being smoothly slice (for example, \cite{Livingstoncomp}). 

Often, we would like to be able to show that a collection of $n$ knots is linearly independent, that is, that they freely generate a subgroup of rank $n$ in $\mathcal{C}$. One way to accomplish this is to define a concordance homomorphism whose domain has rank at least $n$, and to show that the image of this collection of knots has span equal to $n$. Thus, the $\Z$-valued concordance homomorphism $\tau$ is not sufficient for this type of result.

We turn to the more robust invariant $CFK^{\infty}(K)$. In \cite{HomCables}, we defined a $\{-1, 0, 1 \}$-valued concordance invariant, $\varepsilon(K)$. The invariant $\varepsilon$ is associated to the $\Z \oplus \Z$ filtered chain complex $\CFKi$ in a manner similar to how $\tau$ is associated to the $\Z$-filtered chain complex $\widehat{CFK}$; that is, we ask when certain natural maps vanish on homology. We will sometimes write $\varepsilon(\CFKi(K))$, rather than $\varepsilon(K)$, to emphasize that $\varepsilon$ is an invariant associated to the knot Floer complex of $K$.

The goal of this paper is to use $\varepsilon$ to define a new concordance homomorphism that is strong enough to detect linear independence in $\mathcal{C}$.
The main idea is to turn the monoid of chain complexes $CFK^{\infty}(K)$ (under tensor product) into a group, which we will denote $\cF$, in much the same way that the monoid of knots (under connected sum) can be made into the group $\mathcal{C}$ by quotienting by slice knots. 

\begin{Definition}
Let $\CFKi(K)^*$ denote the dual of $\CFKi(K)$. Define the group $\cF$ to be
$$\cF= \big( \{CFK^{\infty}(K) \ | \ K\subset S^3 \}, \otimes \big) / \sim$$
where
$$\CFKi(K_1) \sim \CFKi(K_2) \iff \varepsilon \big(\CFKi(K_1) \otimes \CFKi(K_2)^*\big)=0.$$
\end{Definition}

\begin{Theorem}
The map
$$\mathcal{C} \rightarrow \cF,$$
sending a class in $\mathcal{C}$ represented by $K$ to the class in $\cF$ represented by $CFK^{\infty}(K)$ is a group homomorphism.
\end{Theorem}

This group $\cF$ has the advantage that it can be studied from an algebraic perspective, much like the algebraic concordance group defined by Levine \cite{Levine1, Levine2} in terms of the Seifert form. However, Levine's homomorphism factors through the topological concordance group, while ours does not.

One algebraic feature of $\cF$ is that it is totally ordered, with an additional well-defined notion of domination,``$\ll$''. Moreover, 
we can use the relation $\ll$ to define a filtration on $\cF$ that can be used to show linear independence of certain classes. Given a chain 
$$0<[\CFKi(K_1)] \ll [\CFKi(K_2)] \ll \ldots \ll [\CFKi(K_n)],$$
it follows that the collection
$$ \big\{ [\CFKi(K_i)] \big\}_{i=1}^n $$
is linearly independent in $\cF$, and hence
$$ \big\{ [K_i] \big\}_{i=1}^n $$
is independent in $\cC$. (It is also possible to use spectral sequences to define a second, independent filtration on the group $\cF$.) One consequence of this filtration is that $\cF$ contains a subgroup isomorphic to $\Z ^\infty$; see Theorem \ref{thm:independence} below. We will use this rich structure on $\cF$ to better understand $\cC$.

Let $T_{p,q}$ denote the $(p,q)$-torus knot, $K_{p,q}$ the $(p,q)$-cable of $K$ (where $p$ denotes the longitudinal winding and $q$ denotes the meridional winding), and $D$ the (positive, untwisted) Whitehead double of the right-handed trefoil. We write $T_{m, n; p,q}$ to denote the $(p,q)$-cable of the $(m, n)$-torus knot. Let $-K$ denote the reverse of the mirror image of $K$, that is, the inverse of $K$ in $\cC$.

\begin{Theorem}
\label{thm:independence}
The topologically slice knots
$$D_{p, p+1} \# -T_{p, p+1}, \quad p \geq 1$$
are independent in the smooth concordance group; that is, they freely generate a subgroup of infinite rank.
\end{Theorem}

The first example of an infinite family of smoothly independent, topologically slice knots was given by Endo \cite{Endo}. His examples consist of certain pretzel knots. More recently, Hedden and Kirk \cite{HeddenKirk} showed that an infinite family of (untwisted) Whitehead doubles of certain torus knots are smoothly independent. The structure of $\cF$ shows that our examples (when $p  >1$) are smoothly independent from both of these earlier families.

Let $P(K)$ denote the satellite of $K$ with pattern $P$; that is, $P$ is a knot in $S^1 \times D^2$, which we then glue into the (zero framed) knot complement $S^3 - \textup{nbd } K$ to obtain the knot $P(K) \subset S^3$.
Recall that the map $P(-): \mathcal{C} \rightarrow \mathcal{C}$ given by
$$[K] \mapsto [P(K)]$$
is well-defined, by ``following'' the concordance along the satellite.

We obtain a similar well-defined map on $\cF$:
\begin{Proposition}
\label{prop:satellite}
The map $P(-): \cF \rightarrow \cF$ given by
$$[CFK^{\infty}(K)] \mapsto [CFK^{\infty}\big(P(K)\big)]$$
is well-defined.
\end{Proposition}

By composing $P$ with $\tau$, we obtain a new concordance invariant
$$\tau_P (K) = \tau \big( P(K) \big),$$
since $K_1$ being concordant to $K_2$ implies that $P(K_1)$ is concordant to $P(K_2)$. In the following theorem, we relate this to $[\CFKi(K)]$.

\begin{Theorem}
\label{thm:patterns}
$[\CFKi(K_1)]=[\CFKi(K_2)]$ if and only if $\tau_P(K_1)=\tau_P(K_2)$ for all patterns $P \subset S^1 \times D^2$.
\end{Theorem}

\noindent Recall that $\tau$ is associated to the weaker, $\Z$-filtered chain complex $\widehat{CFK}$. The above theorem says that knowing information about a weaker invariant, namely $\tau$, of satellites of $K$ tells us information about the stronger invariant, $\CFKi$, of the knot itself.

Does the map $P(-): \cC \rightarrow \cC$ always take linearly independent collections of knots to linearly independent collections of knots? We address this question for cables in the following theorem:

\begin{Theorem}
\label{thm:independentsatellites}
For each $n \in \N$, there exists a collection of linearly independent knots
$$\{ K^i \}_{i=1}^n$$
such that for $m \geq n^2-n-1$,
$$\{ K^i_{ 2, 2m+1} \}_{i=1}^n$$ 
is a collection of linearly independent knots in $\cC$.
\end{Theorem}

\noindent This result should be compared to the work of Hedden and Kirk \cite{HeddenKirk}, where they use instantons to prove that the Whitehead doubles of $(2, 2^n+1)$-torus knots are linearly independent.

Central to the definition of $\cF$ is the concordance invariant $\varepsilon$, which exhibits the following properties:

\begin{itemize}
	\item If $K$ is smoothly slice, then $\varepsilon(K)=0$.
	\vspace{4pt}
	\item If $\varepsilon(K)=0$, then $\tau(K)=0$.
	\vspace{4pt}
	\item There exist knots $K$ with $\tau(K)=0$ but $\varepsilon(K) \neq 0$; that is, $\varepsilon$ is \emph{strictly stronger} than $\tau$ at obstructing sliceness.
	\vspace{4pt}
	\item $\varepsilon(-K)=-\varepsilon(K)$.
	\vspace{4pt}
	\item If $\varepsilon(K)=\varepsilon(K')$, then $\varepsilon(K \# K')=\varepsilon(K)$. If $\varepsilon(K)=0$, then $\varepsilon(K \# K')=\varepsilon(K')$.
\end{itemize}
These facts are proved in \cite{HomCables}; we give sketches of their proofs in Section \ref{sec:epsilon}. Notice that since $\varepsilon(K)=0$ implies that $\tau(K)=0$, the map
$$\tau: \cC \rightarrow \Z$$
factors through $\cF$.

\vspace{.5cm}
\noindent \textbf{Organization.} We begin by recounting the necessary definitions and properties of the complex $CFK^{\infty}$ (Section \ref{sec:CFK}) and the concordance invariant $\varepsilon$ (Section \ref{sec:epsilon}). With these definitions in place, we proceed to define the group $\cF$, describe its various algebraic properties, and give examples (Section \ref{sec:F}). We study satellites in Section \ref{sec:satellites}. We conclude with the algebraic details in Section \ref{sec:tedious}.

We work with coefficients in $\F=\Z/2\Z$ throughout.

\vspace{.5cm}
\noindent \textbf{Acknowledgements.} I would like to thank Paul Melvin, Chuck Livingston, Matt Hedden, Rumen Zarev, Robert Lipshitz, Peter Ozsv\'ath, and Dylan Thurston for helpful conversations, and Peter Horn for his comments on an earlier version of this paper.

\section{The knot Floer complex $CFK^{\infty}$}
\label{sec:CFK}

To a knot $K \subset S^3$, Ozsv\'ath and Szab\'o \cite{OSknots}, and independently Rasmussen \cite{R}, associate $CFK^{\infty}(K)$, a $\Z$-filtered chain complex over $\F[U, U^{-1}]$, whose filtered chain homotopy type is an invariant of $K$. The complex $CFK^{\infty}$ can be considered as a $\Z \oplus \Z$-filtered chain complex, with the second filtration induced by $-(U\textup{-exponent})$. The ordering on $\Z \oplus \Z$ is given by $(i, j) \leq (i', j')$ if $i \leq i'$ and $j \leq j'$. We assume the reader is familiar with this invariant, and the various related flavors, $CFK^-(K)$ and $\widehat{CFK}(K)$; for an expository introduction to these invariants, see \cite{OSsurvey}.
The knot $K$ is specified by a doubly pointed Heegaard diagram, $(\Sigma, \boldsymbol{\alpha}, \boldsymbol{\beta}, w, z)$, and the generators (over $\F[U, U^{-1}]$) of $CFK^{\infty}(K)$ are the usual $g$-tuples of intersection points between the $\alpha$- and $\beta$-circles, where $g$ is the genus of $\Sigma$ and each $\alpha$-circle and each $\beta$-circle is used exactly once. The differential is defined as
$$\partial x = \sum_{y \in \mathfrak{S}(\mathcal{H})} \sum_{\substack{\phi \in \pi_2(x, y) \\ \mathrm{ind}(\phi)=1}} \#\widehat{\mathcal{M}}(\phi) \  U^{n_w(\phi)} \cdot y.$$

This complex is endowed with a homological $\Z$-grading, called the \emph{Maslov grading M}, as well as a $\Z$-filtration, called the \emph{Alexander filtration A}. The relative Maslov and Alexander gradings are defined as
$$M(x)-M(y)=\mathrm{ind}(\phi)-2n_w(\phi) \qquad \mathrm{and}  \qquad A(x)-A(y)=n_z(\phi)-n_w(\phi),$$
for $\phi \in \pi_2(x, y)$.
The differential, $\partial$, decreases the Maslov grading by one, and respects the Alexander filtration; that is,
$$M(\partial x)=M(x)-1\qquad \mathrm{and} \qquad A(\partial x) \leq A(x).$$
Multiplication by $U$ shifts the Maslov grading and respects the Alexander filtration as follows:
$$M(U \cdot x)=M(x)-2\qquad \mathrm{and} \qquad A(U \cdot x) = A(x)-1.$$

It is often convenient to view this complex in the $(i, j)$-plane, where the $i$-axis represents $-(U\textup{-exponent})$ and the $j$-axis represents the Alexander filtration. The Maslov grading is suppressed from this picture. We place a generator $x$ at position $(0, A(x))$; more generally, an element of the form $U^i \cdot x$ will have coordinates $(-i, A(x)-i)$.

A basis $\{x_i\}$ for a filtered chain complex $(C, \partial)$ is called a \emph{filtered basis} if the set $\{x_i \ | \ x_i \in C_S \}$ is a basis for $C_S$ for all filtered subcomplexes $C_S \subset C$. Given a filtered basis for $\CFKi$, we may visualize the differential by placing an arrow from a generator $x$ to a generator $y$ if $y$ appears in $\partial x$. The differential points non-strictly to the left and down. Often, it will be convenient to consider only the part of the differential that preserves the Alexander grading, i.e., the horizontal arrows. We will denote this by $\partial^\horz$. Similarly, we will use $\partial^\vert$ to denote the part of the differential that preserves the filtration by powers of $U$, i.e., the vertical arrows.

Given $S \subset \Z \oplus \Z$, let $C\{S\}$ denote the set of elements in the plane whose $(i,j)$-coordinates are in $S$ together with the arrows between them. The complex $CFK^-(K)$ is the subcomplex $C\{i \leq 0\}$, that is, the left half-plane.  The complex $\widehat{CFK}(K)$ is the subquotient complex $C\{i=0\}$. 

The integer-valued smooth concordance invariant $\tau(K)$ is defined in \cite{OS4ball} to be
$$\tau(K) = \mathrm{min} \{ s \ | \ \iota: C\{ i=0, j \leq s \} \rightarrow C\{ i=0 \} \textup{ induces a non-trivial map on homology} \},$$
where $\iota$ is the natural inclusion of chain complexes. Alternatively, $\tau(K)$ may be defined in terms of the $U$-action on $HFK^-(K)$, as in \cite[Appendix A]{OST}:
$$\tau(K)= -\textup{max} \{ s \ | \  \exists \ [x] \in HFK^-(K, s) \textup{ such that } \forall \ d \geq 0, \ U^d[x]\neq 0 \},$$
where $HFK^-(K, s) = H_*(C\{ i \leq 0, j=s \})$.

The complex $CFK^{\infty}(K)$ satisfies certain symmetry and rank properties \cite[Section 3]{OSknots}. The complex obtained by interchanging the roles of $i$ and $j$ is filtered chain homotopic to the original. Also, the rank of the homology of any column or row is one; more generally, modulo grading shifts, any column or row is filtered chain homotopic to $\widehat{CFK}(K)$.

By \cite[Theorem 7.1]{OSknots}, we have the filtered chain homotopy equivalence
$$\CFKi(K_1 \# K_2) \simeq \CFKi(K_1) \otimes_{\F[U, U^{-1}]} \CFKi(K_2).$$
Let $-K$ denote the reverse of the mirror image of $K$. The knot Floer complex is not sensitive to changes in orientation of the knot, but it is sensitive to changes in the orientation of the ambient manifold \cite[Section 3.5]{OSknots}. In particular,
$$\CFKi(-K) \simeq \CFKi(K)^*,$$
where $\CFKi(K)^*$ denotes the dual of $\CFKi(K)$, i.e., $\textup{Hom}_{\F[U, U^{-1}]}(\CFKi(K), \F[U, U^{-1}])$. To depict the complex $\CFKi(K)^*$ in the $(i, j)$-plane, we take the complex $\CFKi(K)$ and reverse the direction of all of the arrows as well as the direction of both of the filtrations. (In practice, we can accomplish this by reversing the direction of all of the arrows and then turning our heads upside down.)

We point out that when we write $\CFKi(K)$, we are really denoting an equivalence class of filtered chain complexes. We may always choose as our representative the $E_1$ page of the spectral sequence associated to one of these complexes, that is, the homology of the associated graded object together with the induced differentials. In other words, we may choose our representative to be reduced, in the sense that any differential strictly lowers the filtration (in at least one direction).

\section{The invariant $\varepsilon$}
\label{sec:epsilon}

The invariant $\varepsilon$ can be defined in terms of the (non-)vanishing of certain cobordism maps, which, using the relation between large surgery and knot Floer homology \cite[Theorems 4.1 and 4.4]{OSknots}, has an algebraic interpretation in terms of the filtered chain complex $CFK^{\infty}(K)$.

Let $N$ be a sufficiently large integer. (It turns out that $N>2g(K)$ will suffice; see \cite[Theorem 1.1]{OSinteger} and \cite[Theorem 5.1]{OSknots}.) We consider the map
$$F_s: \widehat{HF}(S^3) \rightarrow \widehat{HF}(S^3_{-N}(K), [s]),$$
induced by the $2$-handle cobordism, $W^4_{-N}$. As usual, $[s]$ denotes the restriction to $S^3_{-N}(K)$ of the Spin$^c$ structure $\mathfrak{s}_s$ over $W^4_{-N}$ with the property that
$$\langle c_1(\mathfrak{s}_s), [\widehat{F}] \rangle +N = 2s,$$
where $|s| \leq \tfrac{N}{2}$ and $\widehat{F}$ denotes the capped off Seifert surface in the four manifold.
We also consider the map
$$G_s: \widehat{HF}(S^3_N(K), [s]) \rightarrow \widehat{HF}(S^3),$$
induced by the $2$-handle cobordism, $-W^4_N$. 

The maps $F_s$ and $G_s$ can be defined algebraically by studying certain natural maps on subquotient complexes of $CFK^{\infty}(K)$, as in \cite{OSknots}. The map $F_s$ is induced by the chain map
$$C\{i=0 \} \rightarrow C\{ \textup{min} (i, j-s)=0\}$$
consisting of quotienting by $C\{ i=0, j <s \}$, followed by inclusion. Similarly, the map $G_s$ is induced by the chain map
$$C\{ \textup{max } (i, j-s)=0\} \rightarrow C\{i=0 \}$$
consisting of quotienting by $C\{i<0, j=s \}$, followed by inclusion.

For ease of notation, we will often write simply $\tau$ for $\tau(K)$ when the meaning is clear from context. Notice that for $s>\tau$, $F_s$ is trivial, since quotienting $C\{i=0\}$ by $C\{ i=0, j<s \}$ will induce the trivial map, as the homology of $C\{ i=0 \}$ is supported in filtration level $\tau$.

For $s< \tau$, $F_s$ is non-trivial, since any generator of $H_*(C\{ i=0\})$ will still be in the kernel, but not the image, of the differential on $C\{ \textup{min} (i, j-s)=0\}$.

The map $F_\tau$ may be trivial or non-trivial, depending on whether the class representing a generator of $H_*(C\{i=0\})$ lies in the image of the differential on $C\{ \textup{min} (i, j-\tau)=0\}$ or not.

The maps $G_{\tau}$ behaves similarly. For $s>\tau$, the map $G_s$ is non-trivial, and for $s<\tau$, $G_s$ is trivial. The map $G_\tau$ will be non-trivial if the class representing a generator of $H_*(C\{i=0\}$ lies in the kernel of the differential on $C\{ \textup{max } (i, j-s)=0\}$, and trivial otherwise.

Because $C\{j=\tau\}$ is a chain complex, and so $\partial^2=0$, it follows that $F_s$ and $G_s$ cannot both be trivial; that is, a class cannot lie in the image but not in the kernel of the differential. (This is made precise in \cite{HomCables}.) Therefore, there are three possibilities for $F_\tau$ and $G_\tau$: either exactly one vanishes, or neither vanishes.

\begin{definition}
The invariant $\varepsilon$ is defined in terms of $F_\tau$ and $G_\tau$ as follows:

\begin{itemize}
	\item $\varepsilon (K)=1$ if and only if $F_\tau$ is trivial (in which case $G_\tau$ is necessarily non-trivial).
	\item $\varepsilon (K)=-1$ if and only if $G_\tau$ is trivial (in which case $F_\tau$ is necessarily non-trivial).
	\item $\varepsilon (K)=0$ if and only if both $F_\tau$ and $G_\tau$ are non-trivial.
\end{itemize}
\end{definition}

Let $[x]$ be a generator of $H_*(C\{i=0\})$, the so-called ``vertical'' homology. In light of the preceding discussion, the definition of $\varepsilon$ corresponds to viewing $[x]$ as a class in the ``horizontal'' complex $C\{j=\tau\}$ as follows:
\begin{itemize}
	\item $\varepsilon(K)=1$ if and only if $[x]$ is in the image of horizontal differential.
	\item $\varepsilon(K)=-1$ if and only if $[x]$ is not in the kernel of the horizontal differential.
	\item $\varepsilon(K)=0$ if and only if $[x]$ is in the kernel but not the image of the horizontal differential.
\end{itemize}
Notice that $\varepsilon$ is an invariant of the filtered chain homotopy type of $CFK^{\infty}$; at times, to emphasize this point, we will write $\varepsilon(CFK^{\infty}(K))$ rather than simply $\varepsilon(K)$. 

This idea of associating numerical invariants to filtered chain complexes is common; for example, to any $\Z$-filtered chain complex whose total homology has rank one, we can define an integer-valued invariant that measures the minimum filtration level at which this homology is supported, e.g., $\tau(K)$, which is an invariant of the $\Z$-filtered chain homotopy type of $\widehat{CFK}(K)$.

Similarly, to any $\Z \oplus \Z$-filtered chain complex whose ``vertical'' homology has rank one, we can define a $\{-1, 0, 1\}$-valued invariant that measures how this class appears in the ``horizontal'' complex, i.e., in the image of the horizontal differential, in the kernel but not the image, or not in the kernel, respectively. In particular, when $\varepsilon(K)=0$, then $\CFKi(K)$ is filtered chain homotopic to a complex with a distinguished generator that is non-trivial in both the vertical and the horizontal homology.

\begin{proposition}[\cite{HomCables}] \label{prop:propepsilon} The following are properties of $\varepsilon(K)$:
\begin{enumerate}
\item \label{item:slice} If $K$ is smoothly slice, then $\varepsilon(K)=0$.
\vspace{5pt}
\item \label{item:tau0} If  $\varepsilon(K)=0$, then $\tau(K)=0$.
\vspace{5pt}
\item \label{item:mirror} $\varepsilon(-K)=-\varepsilon(K)$.
\vspace{5pt}
\item \label{item:connectedsum}
	\begin{enumerate}
		\item If $\varepsilon(K)=\varepsilon(K')$, then $\varepsilon(K \# K')=\varepsilon(K)=\varepsilon(K')$.
		\vspace{5pt}
		\item If $\varepsilon(K)=0$, then $\varepsilon(K \# K')=\varepsilon(K')$.
	\end{enumerate}
\end{enumerate}
\end{proposition}

\noindent For completeness, we sketch the proof below.

\begin{proof}[Sketch of proof]
To prove (\ref{item:slice}), we consider the $d$-invariants of large surgery along $K$.
If $K$ is slice, then the surgery correction terms defined in \cite{OSabsgr} vanish, i.e., agree with the surgery correction terms of the unknot, and the maps
$$\widehat{HF}(S^3_{N}(K), [0]) \rightarrow \widehat{HF}(S^3) \quad \textrm{and} \quad \widehat{HF}(S^3) \rightarrow \widehat{HF}(S^3_{-N}(K), [0])$$
are non-trivial. Indeed, the surgery corrections terms can be defined in terms of the maps
$${HF}^+(S^3) \rightarrow {HF}^+(S^3_{-N}(K), [s])$$
and we have the commutative diagram
$$\begin{CD}
	 \widehat{HF}(S^3) @>F_s>> \widehat{HF}(S^3_{-N}(K), [s]) \\
	@V\iota VV						@V\iota_{s}VV\\
	 {HF}^+(S^3) @>F^+_s>> {HF}^+(S^3_{-N}(K), [s]).
\end{CD}$$
Let $N\gg0$. If the surgery corrections terms vanish (that is, agree with those of the unknot), then $F^+_\tau$ is an injection \cite[Section 2.2]{RGT}, and so the composition $\iota \circ F^+_\tau$ is non-trivial. By commutativity of the diagram, it follows that $F_\tau$ must be non-trivial. A similar diagram in the case of large positive surgery shows that $G_\tau$ must be non-trivial as well. Hence $\varepsilon(K)=0$.

The proof of (\ref{item:tau0}) follows from the fact that if $\varepsilon(K)=0$, then there is a class $x$ in $\CFKi(K)$ which generates both $H_*(C\{i=0\})$ and $H_*(C\{j=0\})$. In the former complex, $x$ has Alexander grading $A(x)$, and in the latter, viewed as a $\Z$-filtered complex, $x$ has filtration level $-A(x)$. Hence $\tau(K)=-\tau(K)=0$.

The proof of (\ref{item:mirror}) follows from the symmetry properties of the knot Floer complex \cite[Section 3.5]{OSknots}; in particular, we have the filtered chain homotopy equivalence $\CFKi(-K) \simeq \CFKi(K)^*$.

To prove the first part of (\ref{item:connectedsum}): if $[x]$ and $[x']$ are generators of $H_*(\widehat{CFK}(K))$ and $H_*(\widehat{CFK}(K'))$, respectively, then $[x\otimes x']$ is a generator of $H_*(\widehat{CFK}(K \# K'))$. (Here, we are identifying  $\widehat{CFK}$ with $C\{i=0\}$.) Suppose $\varepsilon(K)=\varepsilon(K')=1$. Then both $[x]$ and $[x']$ are both in the image of the horizontal differential, and hence $[x\otimes x']$ is also. The other cases follow similarly.
\end{proof}

\noindent Notice that Proposition \ref{prop:propepsilon} implies that $\varepsilon$ is a concordance invariant. If $K$ and $K'$ are concordant, then $\varepsilon(K \# -K')=0$, in which case $\varepsilon(K)=-\varepsilon(-K')$ by (\ref{item:connectedsum}), or $\varepsilon(K)=\varepsilon(K')$.

Note that we have the following subgroup of $\mathcal{C}$:
$$\{[K] \ | \ \varepsilon(K)=0 \} \subset \mathcal{C}.$$
This observation will useful in the next section.

\section{The group $\cF$}
\label{sec:F}

In this section, we define the group $\cF$ as well as some of its algebraic structure. We will give examples of knots that demonstrate the richness of this structure. In particular, we give an infinite family of topologically slice knots that are linearly independent in $\cF$, and hence also in the smooth concordance group $\cC$, as needed for the proof of Theorem \ref{thm:independence}.
\subsection{Definition of the group $\cF$}
We define the group $\cF$ as
$$\cF = \big( \{ CFK^{\infty}(K) \ | \ K \subset S^3 \}, \otimes \big) / \sim,$$
where
$$CFK^{\infty}(K_1) \sim CFK^{\infty}(K_2) \iff \varepsilon \big(CFK^{\infty}(K_1) \otimes CFK^{\infty}(K_2)^* \big)=0,$$ 
$CFK^{\infty}(K)^*$ denotes the dual of $CFK^{\infty}(K)$, and the tensor product is over $\F[U, U^{-1}]$.
We have the well-defined group homomorphism
$$\cC \rightarrow \cF,$$
given by
$$[K] \mapsto [CFK^{\infty}(K)].$$
Well-definedness follows from the following facts (the first two from \cite[Section 3.5]{OSknots} and the last from Proposition \ref{prop:propepsilon}):
\begin{itemize}
	\item $CFK^{\infty}(-K) \simeq CFK^{\infty}(K)^*$.
	\item $CFK^{\infty}(K_1 \# K_2) \simeq CFK^{\infty}(K_1) \otimes CFK^{\infty}(K_2)$.
	\item If $K$ is smoothly slice, then $\varepsilon(CFK^{\infty}(K))=0$.
\end{itemize}
Notice that $\cF$ is isomorphic to the quotient
$$\cF \cong \cC / \{[K] \ | \ \varepsilon(K)=0 \}.$$
For ease of notation, from now on, we will write
$$\llbracket K \rrbracket$$
to denote $[\CFKi(K)]$, and, when convenient, we will write
$$\llbracket K_1\rrbracket  +\llbracket K_2\rrbracket $$
to denote the operation on the group, which can be thought of as either $[\CFKi(K_1) \otimes \CFKi(K_2)]$ or $[\CFKi(K_1 \# K_2)]$. Note that $-\llbracket K\rrbracket =\llbracket -K\rrbracket $. We denote the identity of $\cF$, $\llbracket \textup{unknot} \rrbracket$, by $0$.

The group $\cF$ has a rich algebraic structure: it has a total ordering, and a ``$\ll$'' relation that satisfies the expected properties and induces a filtration on the group. This algebraic structure on $\cF$ will in turn be useful in understanding the structure of the smooth concordance group $\cC$.

\begin{proposition}
The group $\cF$ is totally ordered, with the ordering given by
$$\llbracket K_1\rrbracket  > \llbracket K_2\rrbracket  \iff \varepsilon(K_1 \# -K_2)=1.$$
\end{proposition}

\begin{proof}
We may think of $\varepsilon(K)$ as the ``sign'' of $\llbracket K\rrbracket $, and then the order relation between any two classes is determined by the sign of their difference.

This relation is clearly transitive, since given
$$\llbracket K_1\rrbracket  > \llbracket K_2\rrbracket  \quad \textup{and} \quad \llbracket K_2\rrbracket  > \llbracket K_3\rrbracket ,$$
it follows that
$$\llbracket K_1\rrbracket  > \llbracket K_3\rrbracket .$$
Indeed,
\begin{align*}
\varepsilon(K_1 \# -K_3)&=\varepsilon(K_1 \# - K_2 \# K_2 \# -K_3) \\
&=1,
\end{align*}
by (\ref{item:connectedsum}) of Proposition \ref{prop:propepsilon} since $\varepsilon(K_1 \# -K_2)=1$ and $\varepsilon(K_2 \# -K_3)=1$.

This relation is also translation invariant. Given
$$\llbracket K_1\rrbracket  > \llbracket K_2\rrbracket ,$$
it follows that
$$\llbracket K_1\rrbracket +\llbracket K_3\rrbracket  > \llbracket K_2\rrbracket +\llbracket K_3\rrbracket ,$$
since
\begin{align*}
\varepsilon(K_1 \# K_3 \# -K_3 \# -K_2) &= \varepsilon(K_1 \# -K_2) \\
&= 1.
\end{align*}
\end{proof}

Totally ordered groups give rise to many natural algebraic constructions, which we will utilize below.
For example, we have a notion of absolute value; that is, given an element $\llbracket K \rrbracket$, either $\llbracket K\rrbracket $ or $-\llbracket K\rrbracket $ is greater than the identity, so we define the absolute value as

$$\big|\llbracket K\rrbracket \big|=\left\{
	\begin{array}{ll}
		\llbracket K\rrbracket  & \text{if } \varepsilon(K) \geq 0 \\
		-\llbracket K\rrbracket  & \text{otherwise}.
	\end{array} \right.$$

A natural question to ask is: Do there exist knots $K_1$ and $K_2$ with $\varepsilon(K_1)=\varepsilon(K_2)=1$ (i.e., they are both ``positive'' with respect to the ordering), and
$$\llbracket K_1]>n\llbracket K_2\rrbracket  \ \ \textup{for all } n \in \N?$$
The answer, it turns out, is yes, motivating the following definition:

\begin{definition}
The class $\llbracket K_1\rrbracket $  \emph{dominates}  $\llbracket K_2\rrbracket $, denoted
$$\llbracket K_1\rrbracket \gg\llbracket K_2\rrbracket ,$$
if $\llbracket K_1]>n\llbracket K_2\rrbracket > 0  \ \ \textup{for all } n \in \N$.
\end{definition}

\noindent Transitivity of $\gg$ follows exactly as for the total ordering. We have the following lemma, showing that the $\gg$ relation satisfies a property we would expect of a ``much bigger'' relation:

\begin{lemma}
\label{lem:gg}
If
$$\llbracket K_1\rrbracket  \gg \llbracket K_2\rrbracket  \quad \textup{ and } \quad \llbracket K_1\rrbracket  \gg \llbracket K_3\rrbracket $$
then
$$\llbracket K_1\rrbracket  \gg \llbracket K_2\rrbracket  +\llbracket K_3\rrbracket .$$
\end{lemma}

\begin{proof}
To see that this is true, we proceed by contradiction. Assume there exists $n \in \N$ such that
$$\llbracket K_1\rrbracket \leq n\big(\llbracket K_2\rrbracket +\llbracket  K_3\rrbracket \big).$$
Then $2\llbracket K_1\rrbracket \leq 2n\big(\llbracket K_2\rrbracket +\llbracket  K_3\rrbracket \big)$, i.e.,
$$\llbracket K_1\rrbracket -2n\llbracket K_2\rrbracket  + \llbracket K_1\rrbracket  -2n\llbracket K_3\rrbracket  \leq 0.$$
But $\llbracket K_1\rrbracket -m\llbracket K_2\rrbracket >0$ and $\llbracket K_1\rrbracket -m\llbracket K_3\rrbracket >0$ for all $m \in \N$, giving us the desired contradiction.
\end{proof}

\begin{remark}\emph{
These ideas could alternatively be phrased in terms of Archimedean equivalence classes. Recall that two elements $a$ and $b$ of a totally ordered group are \emph{Archimedean equivalent} if there exist natural numbers $M$ and $N$ such that $M \cdot |a| > |b|$ and $N \cdot |b| > |a|$. Then we say that $a \gg b$ if $a>b>0$, and $a$ and $b$ are not Archimedean equivalent. Note that the set of Archimedean equivalence classes is naturally totally ordered, and this ordering corresponds  to the $\gg$ relation.
}
\end{remark}

\begin{definition} Let $\cF_K$ denote the collection of elements
$$\cF_K=\{ \llbracket J\rrbracket  \ \big| \ | \llbracket J\rrbracket | \ll | \llbracket K\rrbracket  | \}.$$
\end{definition}

\begin{proposition}
$\cF_K$ is a subgroup of $\cF$.
\end{proposition}

\begin{proof}
If $\llbracket J\rrbracket $ is in $\cF_K$, then $-\llbracket J\rrbracket $ clearly is as well.
Given $\llbracket J_1\rrbracket $ and $\llbracket J_2\rrbracket $ in $\cF_K$, is follows immediately that $\llbracket J_1\rrbracket +\llbracket J_2\rrbracket $ is also in $\cF_K$, by Lemma \ref{lem:gg}.
\end{proof}

\noindent Notice that given a sequence of knots $K_1, K_2, \ldots, K_n$ satisfying
$$\llbracket K_1\rrbracket  \gg \llbracket K_2\rrbracket  \gg \ldots \gg \llbracket K_n\rrbracket ,$$
we obtain a filtration
$$\cF_{K_1} \supset \cF_{K_2} \supset \ldots \supset \cF_{K_n}.$$

\begin{lemma}
\label{lem:ggindependence}
If $\llbracket K_1\rrbracket  \gg \llbracket K_2\rrbracket  \gg \ldots \gg \llbracket K_n\rrbracket > 0$,
then the knots
$$K_1, \ K_2, \ \ldots, K_n$$
are linearly independent in $\cF$ and hence in $\cC$; that is, they generate a subgroup of rank $n$ in both $\cF$ and $\cC$.
\end{lemma}

\begin{proof}
By Lemma \ref{lem:gg}, for any positive integer $m$, $m\llbracket K_1\rrbracket $ dominates any linear combination of $\llbracket K_2\rrbracket , \ldots, \llbracket K_n\rrbracket $, and thus cannot be expressed as a linear combination of these classes. Similarly, $m\llbracket K_i\rrbracket $ dominates any linear combination of $\llbracket K_{i+1}\rrbracket , \ldots, \llbracket K_n\rrbracket $, for $i<n$.
\end{proof}

\subsection{Examples}
We now give examples of families of knots that can be shown to independent in $\cC$. 
\begin{proposition}
\label{prop:examples}
Let $0<p<q$. Then we have the following relations in the group $\cF$:
\begin{enumerate}
	\item $\llbracket T_{p, p+1}\rrbracket  \ll \llbracket T_{q, q+1}\rrbracket $
	\item \label{it:cableofdoublepq} $\llbracket D_{p, p+1}\rrbracket  \ll \llbracket D_{q, q+1}\rrbracket $
	\item \label{it:torusdoublepp} $\llbracket T_{p, p+1}\rrbracket  \ll \llbracket D_{p, p+1}\rrbracket $
	\item $\llbracket T_{p, p+1; 2, 2m+1}\rrbracket  \ll \llbracket T_{q, q+1; 2, 2m+1}\rrbracket $, $m \geq q^2-q-1$.
\end{enumerate}
\end{proposition}
\noindent We will prove this proposition at the end of Section \ref{sec:tedious}.

\begin{remark}
\label{rem:examples}
\emph{A straightforward consequence of (\ref{it:cableofdoublepq}) and (\ref{it:torusdoublepp}) of the preceding proposition is the relation
$$\llbracket D_{p, p+1} \# -T_{p, p+1}\rrbracket  \ll \llbracket D_{q, q+1} \# -T_{q, q+1}\rrbracket .$$
}
\end{remark}

We are now ready to prove Theorem \ref{thm:independence}; that is, we will show that the knots
$$D_{p, p+1} \# -T_{p, p+1}, \quad p \geq 1$$
are smoothly independent while being topologically slice.

\begin{proof}[Proof of Theorem \ref{thm:independence}]
Recall that $D$ is the (positive, untwisted) Whitehead double of the right-handed trefoil. The Alexander polynomial of $D$ is equal to one, and so by Freedman \cite{Freedman}, it follows that $D$ is topologically slice. Hence, the $(p, p+1)$-cable of $D$, $D_{p, p+1}$, is topologically concordant to the $(p, p+1)$-cable of the unknot, i.e., the torus knot $T_{p, p+1}$. Thus, $D_{p, p+1} \# -T_{p, p+1}$ is topologically slice.

It follows from Lemma \ref{lem:ggindependence} and Remark  \ref{rem:examples} that the knots
$$D_{p, p+1} \# -T_{p, p+1}, \quad p \geq 1$$
are linearly independent in $\cF$, and hence also in $\cC$.
\end{proof}

\begin{proof}[Proof of Theorem \ref{thm:independentsatellites}]
We need to find a collection of linearly independent knots $\{ K^i \}_{i=1}^{n}$ such that the collection $\{ K^i_{2, 2m+1} \}_{i=1}^{n}$ is also linearly independent for sufficiently large $m$.

Let $K^i = T_{i, i+1}$, and consider the $(2, 2m+1)$-cable of $K^i$, where $m \geq n^2-n-1$. By Lemma \ref{lem:ggindependence} and Proposition \ref{prop:examples}, it follows that the collection
$$\{ K^i \}_{i=1}^n$$
is linearly independent in $\cF$, hence also in $\cC$. Again, by Lemma \ref{lem:ggindependence} and Proposition \ref{prop:examples}, the collection
$$\{ K^i_{2, 2m+1} \}_{i=1}^n,$$
is also linearly independent in $\cF$ and thus in $\cC$.
\end{proof}

\section{Satellites and $\cF$}
\label{sec:satellites}

Recall that $P(K)$ denotes the satellite of $K$ with pattern $P$; that is, $P$ is a knot in $S^1 \times D^2$, which we then glue into the (zero framed) knot complement $S^3 - \textup{nbd } K$. The map $P(-): \mathcal{C} \rightarrow \mathcal{C}$ given by
$$[K] \mapsto [P(K)]$$
is well-defined, by ``following'' the concordance along the satellite. We will show that an analogous result holds for the group $\cF$.

\begin{proposition}
\label{prop:satellitemap}
The map $P(-): \cF \rightarrow \cF$ given by
$$\llbracket K\rrbracket  \mapsto \llbracket P(K)\rrbracket $$
is well-defined.
\end{proposition}

The following theorem from \cite{HomCables} gives a formula for $\tau(K_{p,q})$ in terms of $\tau(K)$, $\varepsilon(K)$, $p$, and $q$:

\begin{theorem} [\cite{HomCables}]
\label{thm:cables}
Let $K \subset S^3$, and let $p$, $q$ be relatively prime integers with $p>0$. Then the behavior of $\tau(K_{p,q})$ is completely determined by $p$, $q$, $\tau(K)$, and $\varepsilon(K)$. More precisely:
\begin{enumerate}
	 \item If $\varepsilon(K)=1$, then $\tau(K_{p,q})=p\tau(K)+\frac{(p-1)(q-1)}{2}$.
	 \item If $\varepsilon(K)=-1$, then $\tau(K_{p,q})=p\tau(K)+\frac{(p-1)(q+1)}{2}$.
	\item If $\varepsilon(K)=0$, then $\tau(K)=0$ and $\tau(K_{p, q})= \tau(T_{p,q})=\left\{
	\begin{array}{ll}
		\frac{(p-1)(q+1)}{2} & \text{if } q<0\\
		\frac{(p-1)(q-1)}{2} & \text{if } q>0.
	\end{array} \right.$
\end{enumerate}
\end{theorem}

\noindent We see that knowing $\tau(K_{2, 1})$ and $\tau(K_{2, -1})$ is sufficient to determine $\varepsilon(K)$. More precisely,
\begin{itemize}
	\item If $\tau(K_{2, 1})$ is odd, then $\varepsilon(K)=-1$.
	\item If $\tau(K_{2, -1})$ is odd, then $\varepsilon(K)=1$.
	\item If $\tau(K_{2, 1})=\tau(K_{2, -1})=0$, then $\varepsilon(K)=0$.
\end{itemize}
The proof of Proposition \ref{prop:satellitemap} will rely on this observation.

The proof will also rely on facts from bordered Heegaard Floer homology, as defined by Lipshitz, Ozsv\'ath and Thurston \cite{LOT}. We will need only a special case of the formal properties of these invariants, which we recount here.

To a framed knot complement $Y_K$, we associate a left differential graded module $\widehat{CFD}(Y_K)$, whose homotopy equivalence class is an invariant of the framed knot complement \cite[Theorem 1.1]{LOT}. Furthermore, the homotopy equivalence class is completely determined by the complex $CFK^{\infty}(K)$ and the framing $n$ \cite[Theorem 11.27 and A.11]{LOT}. For our purposes here, it will be sufficient to let $Y_K$ be the zero framed knot complement. In \cite{HomCables}, it is shown that if $\varepsilon(K)=0$, then 
$$\widehat{CFD}(Y_{J \# K}) \simeq \widehat{CFD}(Y_J) \oplus A,$$
 for some left differential graded module $A$.

To a knot $P$ in $S^1 \times D^2$, we associate a right $\mathcal{A}_\infty$-module $CFA^-(S^1 \times D^2, P)$. Let $gCFK^-(K)$ denote the associated graded complex of $CFK^-(K)$, i.e.,  $\oplus_s C \{ i \leq 0, j=s \}$. Notice that $HFK^-(K) \cong H_*(gCFK^-(K))$. Then the pairing theorem for bordered Heegaard Floer homology \cite[Theorem 11.21]{LOT} states that we have the following graded chain homotopy equivalence:
$$gCFK^-(S^3, P(K)) \simeq CFA^-(S^1 \times D^2, P) \widetilde{\otimes} \widehat{CFD}(Y_K),$$
where we choose the zero framing for the knot complement $Y_K$, and where $\widetilde{\otimes}$ denotes the $\mathcal{A}_\infty$-tensor product, a generalization of the derived tensor product. In particular, $\widetilde{\otimes}$ respects summands.

\begin{proof}[Proof of Proposition \ref{prop:satellitemap}]
Assume $\varepsilon(K \# -J)=0$. We would like to show that 
$$\varepsilon(P(K) \# -P(J))=0.$$ 
Utilizing the observation above, it is sufficient to show that 
$$\tau \big((P(K) \# -P(J))_{2, \pm 1}\big)=0.$$

\begin{figure}[htb!]
\labellist
\small \hair 2pt
\endlabellist
\centering
\vspace{10pt}
\includegraphics[scale=1.2]{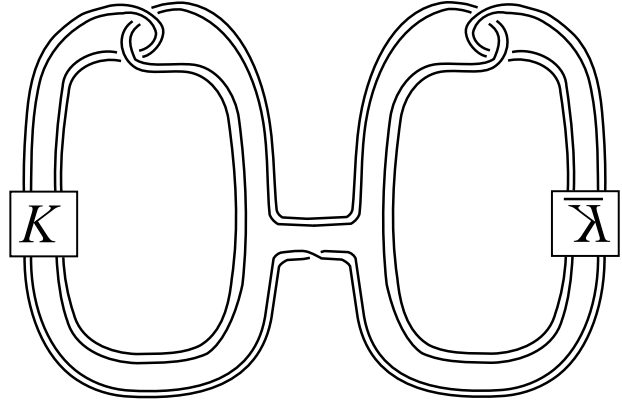}
\vspace{10pt}
\caption{The knot $\big(P(K) \# -P(K)\big)_{2, 1}$, in the case where $P$ is the pattern for the Whitehead double.}
\label{fig:P(K)}
\end{figure}

\begin{figure}[htb!]
\labellist
\small \hair 2pt
\endlabellist
\centering
\vspace{10pt}
\includegraphics[scale=1.2]{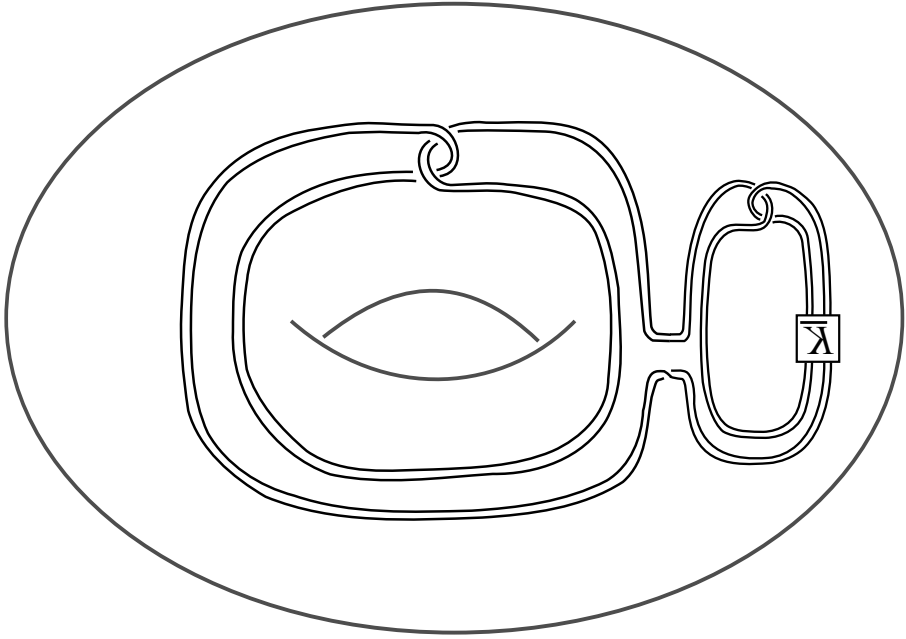}
\vspace{10pt}
\caption{The knot $\big(P(U)\# -P(K)\big)_{2,  1}$ embedded in $S^1 \times D^2$ as the pattern knot $Q$, where again, $P$ is the pattern for the Whitehead double.}
\label{fig:Q(K)}
\end{figure}

Let $U$ denote the unknot. There exists an embedding $Q$ of $\big(P(U)\# -P(J)\big)_{2, \pm 1}$ into $S^1 \times D^2$ such that 
$$Q(K)=\big(P(K) \# -P(J)\big)_{2, \pm 1}.$$
See Figure \ref{fig:Q(K)}. We consider the bordered invariant
$$CFA^-(S^1 \times D^2, Q)$$
associated to $(S^1 \times D^2, Q)$. Notice that $Q(J)=\big(P(J) \# -P(J)\big)_{2, \pm 1}$ is slice and so $\tau\big((P(J) \# -P(J))_{2, \pm 1}\big)=0$.

The knot $K$ is concordant to the knot $K'=J\# K \# -J$. Since $\varepsilon(K \# -J)=0$, we have the following chain homotopy equivalence:
$$\widehat{CFD}(Y_{K'})\simeq \widehat{CFD}(Y_{J}) \oplus A,$$
for some $A$.

The knot $Q(K)$ is concordant to $Q(K')$, since $K$ is concordant to $K'$. The invariant $\tau(Q(K'))$ is determined by 
\begin{align*}
gCFK^-(Q(K')) &\simeq CFA^-(S^1 \times D^2, Q) \widetilde{\otimes} \widehat{CFD}(Y_{K'})\\
&\simeq CFA^-(S^1 \times D^2, Q) \widetilde{\otimes}  \big(\widehat{CFD}(Y_{J}) \oplus A \big) \\
&\simeq gCFK^-\big(Q(J)\big) \oplus B
\end{align*}
where $B$ is the complex $CFA^-(S^1 \times D^2, Q) \widetilde{\otimes}  A$. Notice that $H_*(B)$ is $U$-torsion, since the ranks of $HFK^-(Q(K'))$ and $HFK^-(Q(J))$ as $\F[U]$-modules are both one.
 Thus,
$$\tau(Q(K))=\tau(Q(K'))=\tau(Q(J))=0,$$
since $Q(J)$ is slice. Recalling that $Q(K)=\big( P(K) \# -P(J) \big)_{2, \pm1}$, we have that
$$\tau \big( (P(K) \# -P(J))_{2, \pm 1} \big)=0,$$
implying that
$$\varepsilon(P(K) \# -P(J))=0,$$
as desired.
\end{proof}

We now prove Theorem \ref{thm:patterns}, which we restate here:

\begin{theorem}
$\llbracket K\rrbracket =\llbracket J\rrbracket $ if and only if $\tau_P(K) = \tau_P(J) $ for all patterns $P \subset S^1 \times D^2$.
\end{theorem}

\begin{proof}
The forward direction is true by Proposition \ref{prop:satellitemap} and the fact that the map $\tau: \cC \rightarrow \Z$ factors through $\cF$. We must now show that if $\llbracket K\rrbracket  \neq \llbracket J\rrbracket $, then there exists some pattern $P$ such that $\tau(P(K)) \neq \tau (P(J))$.

Without loss of generality, we may assume that $\varepsilon(K \# -J)=-1$. Let
$$P(K)=(K \# -J)_{2, 1}.$$
Then Theorem \ref{thm:cables} tells us that
$$\tau(P(K))=2(\tau(K) -\tau(J))+1 \qquad \textup{ and } \qquad \tau(P(J))=0,$$
as desired.
\end{proof}

\section{Calculations and a refinement of $\varepsilon$}
\label{sec:tedious}
An element of $\cF$ is an equivalence class of filtered chain complexes. The goal of this section is to define more tractable invariants associated to such a class, compute these invariants for a few families of knots, and show that these invariants are related to the algebraic structure, namely the $\gg$ relation, on $\cF$.

To this end, we will define a refinement of $\varepsilon$. Recall that $\varepsilon$ is defined in terms of whether or not certain maps on subquotient complexes of $\CFKi$ vanish on homology. Our refinement of $\varepsilon$ will be defined in a similar manner.

The invariant $\varepsilon(K)$ is equal to one when the class generating the ``vertical'' homology of $\CFKi(K)$ lies in the image of the horizontal differential. We would like a well-defined way to measure the ``length'' of the differential that hits that class, that is, how much it decreases the horizontal filtration. We will do this by examining certain natural maps on subquotients of $\CFKi$.

The definition of $\varepsilon$ involved examining the map $F_\tau$ induced by
$$C\{i=0\} \rightarrow C\{ \textup{min} (i, j-\tau)=0\}.$$
In particular, if $F_\tau$ is trivial, then $\varepsilon(K)=1$. Consider now the map $H_s$ induced on homology by
$$C\{i=0\} \rightarrow C\{\textup{min} (i, j-\tau)=0, i \leq s\},$$
for some non-negative integer $s$. Notice that $H_0$ is non-trivial, and for sufficiently large $s$, $H_s$ agrees with $F_\tau$.

Suppose that $\varepsilon(K)=1$; that is, $F_\tau$ is trivial. Then define $a_1(K)$ to be
$$a_1(K) =\textup{min} \{ s \ | \ H_s \textup{ is trivial}\}.$$
The idea is that when $\varepsilon(K)=1$, the class generating the vertical homology lies in the image of the horizontal differential, and $a_1$ is measuring the ``length'' of the horizontal differential hitting that class.

Now consider the map $H_{a_1, s}$ induced on homology by
$$C\{i=0\} \rightarrow C\big\{ \{\textup{min} (i, j-\tau)=0, i \leq a_1\} \cup \{ i=a_1, \tau-s \leq j < \tau \} \big\},$$
for some non-negative integer $s$. Clearly, $H_{a_1, 0}$ is trivial. Define
$$a_2(K)=\textup{min} \{ s \ | \ H_{a_1, s} \textup{ is non-trivial}\}.$$
Notice that $a_2(K)$ may be undefined; that is, the map $H_{a_1, s}$ may be trivial for all $s$. Effectively, $a_2$ is measuring the ``length'' of a certain vertical differential, if it exists.

\begin{lemma}
\label{lem:class}
The invariants $a_1$ and $a_2$ are invariants of the class $\llbracket K\rrbracket $.
\end{lemma}

\begin{proof}
Suppose $\llbracket J\rrbracket =\llbracket K\rrbracket $. Then
$$\llbracket J\rrbracket =\llbracket K\rrbracket =\llbracket K \# -J \# J\rrbracket .$$
Since $\varepsilon(K \# -J)=0$, it follows from \cite[Lemma 3.3]{HomCables} that there exists a basis for $\CFKi(K \# -J)$ with a distinguished element, say $x_0$, with no incoming or outgoing horizontal or vertical arrows. Similarly, there is a basis for $\CFKi(J \# -J)$ with a distinguished element $y_0$. Then we may compute $a_1(K \# -J \# J)$ and $a_2(K \# -J \# J)$ by considering either
$$\{x_0 \} \otimes \CFKi(J) \quad \textup{ or } \quad \CFKi(K) \otimes  \{y_0 \},$$
the former giving us $a_1(J)$ and $a_2(J)$, and the latter giving us $a_1(K)$ and $a_2(K)$.
\end{proof}

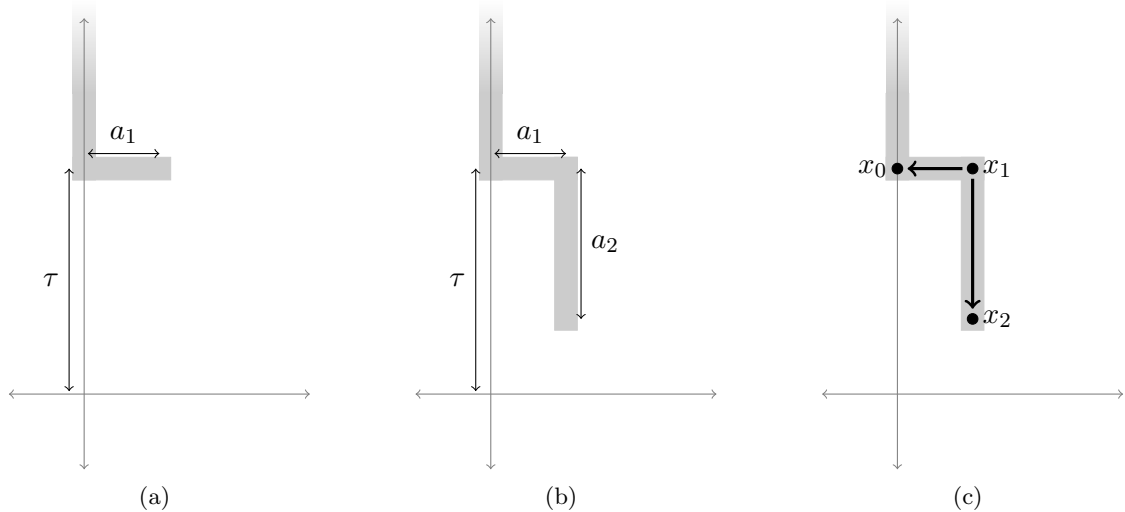
\begin{figure}[htb!]
\vspace{5pt}
\labellist
\small \hair 2pt
\endlabellist
\centering
\subfigure[]{
\begin{tikzpicture}
	\filldraw[black!20!white] (-0.15, 2.85) rectangle (0.15, 4);
	\filldraw[black!20!white] (-0.15, 2.85) rectangle (1.15, 3.15);
	\shade[top color=white, bottom color=black!20!white] (-0.15, 4) rectangle (0.15, 5.25);
	\begin{scope}[thin, gray]
		\draw [<->] (-1, 0) -- (3, 0);
		\draw [<->] (0, -1) -- (0, 5);
	\end{scope}
	\draw [<->] (-0.2, 0.04) -- (-0.2, 3) node [midway, left] {$\tau$};
	\draw [<->] (0.05, 3.2) -- (1, 3.2) node [midway, above] {$a_1$};
\end{tikzpicture}
}
\hspace{25pt}
\subfigure[]{
\begin{tikzpicture}
	\filldraw[black!20!white] (-0.15, 2.85) rectangle (0.15, 4);
	\filldraw[black!20!white] (-0.15, 2.85) rectangle (1.15, 3.15);
	\shade[top color=white, bottom color=black!20!white] (-0.15, 4) rectangle (0.15, 5.25);
	\filldraw[black!20!white] (0.85, 0.85) rectangle (1.15, 3.15);
	\begin{scope}[thin, gray]
		\draw [<->] (-1, 0) -- (3, 0);
		\draw [<->] (0, -1) -- (0, 5);
	\end{scope}
	\draw [<->] (-0.2, 0.04) -- (-0.2, 3) node [midway, left] {$\tau$};
	\draw [<->] (0.05, 3.2) -- (1, 3.2) node [midway, above] {$a_1$};
	\draw [<->] (1.2, 1) -- (1.2, 3) node [midway, right] {$a_2$};	
\end{tikzpicture}
}
\hspace{25pt}
\subfigure[]{
\begin{tikzpicture}
	\filldraw[black!20!white] (-0.15, 2.85) rectangle (0.15, 4);
	\filldraw[black!20!white] (-0.15, 2.85) rectangle (1.15, 3.15);
	\shade[top color=white, bottom color=black!20!white] (-0.15, 4) rectangle (0.15, 5.25);
	\filldraw[black!20!white] (0.85, 0.85) rectangle (1.15, 3.15);
	\begin{scope}[thin, gray]
		\draw [<->] (-1, 0) -- (3, 0);
		\draw [<->] (0, -1) -- (0, 5);
	\end{scope}
	\filldraw (0, 3) circle (2pt) node[] (x_0){};
	\filldraw (1, 3) circle (2pt) node[] (x_1){};
	\filldraw (1, 1) circle (2pt) node[] (x_2){};
	\draw [very thick, <-] (x_0) -- (x_1);
	\draw [very thick, <-] (x_2) -- (x_1);
	\node [left] at (x_0) {$x_0$};
	\node [right] at (x_1) {$x_1$};
	\node [right] at (x_2) {$x_2$};
\end{tikzpicture}
}

\caption{Left, the complex $A$ in the $(i, j)$-plane. Center, the complex $B$. Right, part of the basis in Lemma \ref{lem:basis}.}
\label{fig:basis}
\end{figure}

\begin{lemma}
\label{lem:basis}
Let $a_1=a_1(K)$. Then there exists a basis $\{x_i\}$  over $\F[U, U^{-1}]$ for $\CFKi$ with basis elements $x_0$ and $x_1$ with the property that
\begin{enumerate}
	\item \label{it:a_1} There is a horizontal arrow of length $a_1$ from $x_1$ to $x_0$.	
	\item There are no other horizontal or vertical arrows to or from $x_0$.
	\item There are no other horizontal arrows to or from $x_1$.
	\newcounter{enumi_saved}
	\setcounter{enumi_saved}{\value{enumi}}
\end{enumerate}
If we also have that $a_2=a_2(K)$ is well-defined, then there exists a basis $\{x_i\}$  with basis elements $x_0$, $x_1$, and $x_2$ with the  following properties, in addition to the ones above:
\begin{enumerate}
\setcounter{enumi}{\value{enumi_saved}}
	\item \label{it:a_2} There is a vertical arrow of length $a_2$ from $x_1$ to $x_2$.
	\item There are no other vertical arrows to or from $x_1$ or $x_2$.
\end{enumerate}
\end{lemma}

\begin{proof}
We will give the proof for the case where $a_2$ is well-defined. The proof in the case where $a_2$ is not well-defined is a straightforward simplification of this proof.

For ease of notation, let 
\begin{align*}
A&=C\{\textup{min} (i, j-\tau)=0, i \leq a_1 \}\\
B&=C\big\{ \{\textup{min} (i, j-\tau)=0, i \leq a_1\} \cup \{ i=a_1, \tau-a_2 \leq j < \tau \} \big\},
\end{align*}
so that $H_{a_1}$ and $H_{a_1, a_2}$, respectively, are the maps on homology induced by
\begin{align*}
C\{i=0\} &\rightarrow A\\
C\{i=0\} &\rightarrow B.
\end{align*}
See Figure \ref{fig:basis}. Since $H_{a_1}$ is trivial, it follows that there is a generator, say $x_0$, of $H_*(C\{i=0\})$ in position $(0, \tau)$ that is in the image of the differential on $A$, but not in the image of the differential on $B$. Since $H_{a_1, a_2}$ is non-trivial, there exists a class $x_1$ supported in position $(a_1, \tau)$ whose boundary, in $A$, is $x_0$, and whose boundary, in $B$, is a class, say $x_0+x_2$, where $x_2$ is supported in position $(a_1, \tau-a_2)$. Moreover, we may replace $x_0$ with $\partial^\horz x_1$, since a priori, $\partial^\horz x_1$ might include elements with negative $i$-coordinate. Similarly, we may replace $x_2$ with $\partial^\vert x_1$.

We now complete $\{x_0, x_1, x_2\}$ to a basis $\{x_i\}$ for $\CFKi(K)$, and conditions (\ref{it:a_1}) and (\ref{it:a_2}) above are satisfied. To satisfy the remaining three conditions, we will use a change of basis in order to remove the unwanted arrows.

There are no vertical arrows leaving $x_0$, since it is in the kernel of the vertical differential. Since $x_0$ is not in the image of the vertical differential, if there is an incoming vertical arrow to $x_0$ from, say, $y$, then there is also a vertical arrow from $y$ to, say, $z$. Changing basis to replace $z$ with $z+x_0$ will remove the vertical arrow to $x_0$. All of the incoming vertical arrows to $x_0$ may be removed in this manner, and filtration considerations ensure that we have not changed $x_1$ or $x_2$.

Since $x_0$ is in the image of $\partial^\horz$, it follows immediately that there are no horizontal arrows leaving $x_0$, by the fact that ${\partial^\horz}\circ{\partial^\horz}=0$. We must now remove any horizontal arrows entering $x_0$. Suppose there is an arrow of length $\ell$ from $y$ to $x_0$. If $\ell < a_1$, we may remove the arrow as in the preceding paragraph. If $\ell \geq a_1$, then we replace $y$ with $y+x_1$. In this manner, we can remove all of other horizontal arrows into $x_0$.

There are now no horizontal arrows entering $x_1$, because $\partial^\horz x_1=x_0$, ${\partial^\horz}\circ{\partial^\horz}=0$, and there are no other horizontal arrows to $x_0$.

We may remove unwanted vertical arrows involving $x_1$ and $x_2$ in the same manner that we removed unwanted horizontal arrows involving $x_1$ and $x_0$.
\end{proof}

Note that if we have such a basis $\{x_i\}$ for $\CFKi(K)$, then we have a basis $\{x_i^*\}$ for $\CFKi(K)^*$ satisfying the following:
\begin{itemize}
	\item There is a horizontal arrow of length $a_1(K)$ from $x^*_0$ to $x^*_1$.	
	\item There is a vertical arrow of length $a_2(K)$ from $x^*_2$ to $x^*_1$.
	\item There are no other horizontal or vertical arrows to or from $x^*_0$.
	\item There are no other horizontal or vertical arrows to or from $x^*_1$.
	\item There are no other vertical arrows to or from $x^*_2$.
\end{itemize}
If $x_k$ has filtration level $(i, j)$, then $x_k^*$ has filtration level $(-i, -j)$.
We will use these types of bases to prove the following lemmas:

\begin{lemma}
\label{lem:agg1}
If $a_1(J)>a_1(K)$, then
$$\llbracket K\rrbracket  \gg \llbracket J\rrbracket .$$
\end{lemma}

\begin{proof}
We proceed using induction. We will show that $\varepsilon(K \# -J)=1$ and that
\begin{align*}
a_1(K \# -J) &=a_1(K)
\end{align*}
from which we can conclude that
$$\varepsilon(K \# -nJ)=1$$
for all $n \in \N$.

Let $\{x_i\}$ be a basis for $\CFKi(K)$ found using the first part of Lemma \ref{lem:basis}. Similarly, let $\{y_i\}$ be such a basis for $\CFKi(J)$, and hence $\{y_i^*\}$ is a basis for $\CFKi(-J)$. We consider the knot $K \# -J$ and its knot Floer complex. Notice that $x_0 y_0^* $ generates $H_*(C\{i=0\})$, the ``vertical'' homology of $\CFKi(K \# -J)$. Let $\tau=\tau(K \# -J)$. 

Consider the subquotient complex
$$A=C\{ \textup{min} (i, j-\tau)=0\}.$$
There is a direct summand of $A$ consisting of generators $x_0 y_0^* $ and $x_1 y_0^* $, with a horizontal arrow of length $a_1(K)$ from the latter to the former. Hence, $\varepsilon(K \# -J)=1$ and $a_1(K \# -J) =a_1(K)$, as desired.
\end{proof}

\begin{lemma}
\label{lem:agg}
If $a_1(J)=a_1(K)$ and $a_2(J) > a_2(K)$, then
$$\llbracket J\rrbracket  \gg \llbracket K\rrbracket .$$
\end{lemma}

\begin{figure}[htb!]
\vspace{5pt}
\labellist
\small \hair 2pt
\endlabellist
\centering
\subfigure[]{
\begin{tikzpicture}
	\filldraw[black!20!white] (-0.15, 2.85) rectangle (0.15, 5);
	\filldraw[black!20!white] (-0.15, 2.85) rectangle (1.15, 3.15);
	\shade[top color=white, bottom color=black!20!white] (-0.15, 5) rectangle (0.15, 6.25);
	\filldraw[black!20!white] (0.85, 0.85) rectangle (1.15, 3.15);

	\filldraw (0, 3) circle (2pt) node[] (x_0){};	
	\filldraw (1, 3) circle (2pt) node[] (x_1){};
	\filldraw (1, 1) circle (2pt) node[] (x_2){};
	\draw [very thick, <-] (x_0) -- (x_1);
	\draw [very thick, <-] (x_2) -- (x_1);
	\node [left] at (x_0) {$x_0$};
	\node [right] at (x_1) {$x_1$};
	\node [right] at (x_2) {$x_2$};
	\draw (0, -0.15);
\end{tikzpicture}
}
\hspace{25pt}
\subfigure[]{
\begin{tikzpicture}

	\filldraw (1, 2) circle (2pt) node[] (x_0){};	
	\filldraw (0, 2) circle (2pt) node[] (x_1){};
	\filldraw (0, 4) circle (2pt) node[] (x_2){};
	\draw [very thick, ->] (x_0) -- (x_1);
	\draw [very thick, ->] (x_2) -- (x_1);
	\node [right] at (x_0) {$x^*_0$};
	\node [left] at (x_1) {$x^*_1$};
	\node [right] at (x_2) {$x^*_2$};
	\draw (0, -0.15);
\end{tikzpicture}
}
\hspace{25pt}
\subfigure[]{
\begin{tikzpicture}
	\filldraw[black!20!white] (-0.15, 2.85) rectangle (0.15, 5);
	\filldraw[black!20!white] (-0.15, 2.85) rectangle (1.15, 3.15);
	\shade[top color=white, bottom color=black!20!white] (-0.15, 5) rectangle (0.15, 6.25);
	\filldraw[black!20!white] (0.85, -0.15) rectangle (1.15, 3.15);

	\filldraw (0, 3) circle (2pt) node[] (y_0){};	
	\filldraw (1, 3) circle (2pt) node[] (y_1){};
	\filldraw (1, 0) circle (2pt) node[] (y_2){};
	\draw [very thick, <-] (y_0) -- (y_1);
	\draw [very thick, <-] (y_2) -- (y_1);
	\node [left] at (y_0) {$y_0$};
	\node [right] at (y_1) {$y_1$};
	\node [right] at (y_2) {$y_2$};
\end{tikzpicture}
}
\hspace{25pt}
\subfigure[]{
\begin{tikzpicture}
	\filldraw[black!20!white] (-0.15, 2.85) rectangle (0.15, 5);
	\filldraw[black!20!white] (-0.15, 2.85) rectangle (1.15, 3.15);
	\shade[top color=white, bottom color=black!20!white] (-0.15, 5) rectangle (0.15, 6.25);
	\filldraw[black!20!white] (0.85, -0.15) rectangle (1.15, 3.15);

	\filldraw (0, 3.12) circle (2pt) node[] (x_1y_1){};
	\filldraw (0, 2.88) circle (2pt) node[] (x_0y_0){};	
	\filldraw (1, 3) circle (2pt) node[] (x_0y_1){};
	\filldraw (1, 0) circle (2pt) node[] (x_0y_2){};
	\filldraw (0, 5) circle (2pt) node[] (x_2y_1){};
	\draw [very thick, <-] (x_1y_1) -- (x_0y_1);
	\draw [very thick, <-] (x_0y_0) -- (x_0y_1);
	\draw [very thick, <-] (x_0y_2) -- (x_0y_1);
	\draw [very thick, <-] (x_1y_1) -- (x_2y_1);	
	\node [left] at (x_1y_1) {$x^*_1y_1$};
	\node [left, below] at (x_0y_0) {$x^*_0y_0$};
	\node [right] at (x_0y_1) {$x^*_0y_1$};
	\node [right] at (x_0y_2) {$x^*_0y_2$};
	\node [right] at (x_2y_1) {$x^*_2y_1$};
\end{tikzpicture}
\label{subfig:summand}
}

\caption{Far left, a portion of the basis $\{x_i \}$ for $\CFKi(K)$, followed by a portion of the basis $\{x^*_i \}$ for $\CFKi(K)^*$. Next, a portion of the basis $\{y_i \}$ for $\CFKi(J)$. Far right, a direct summand of the subquotient complex $B$.}
\label{fig:bases}
\end{figure}
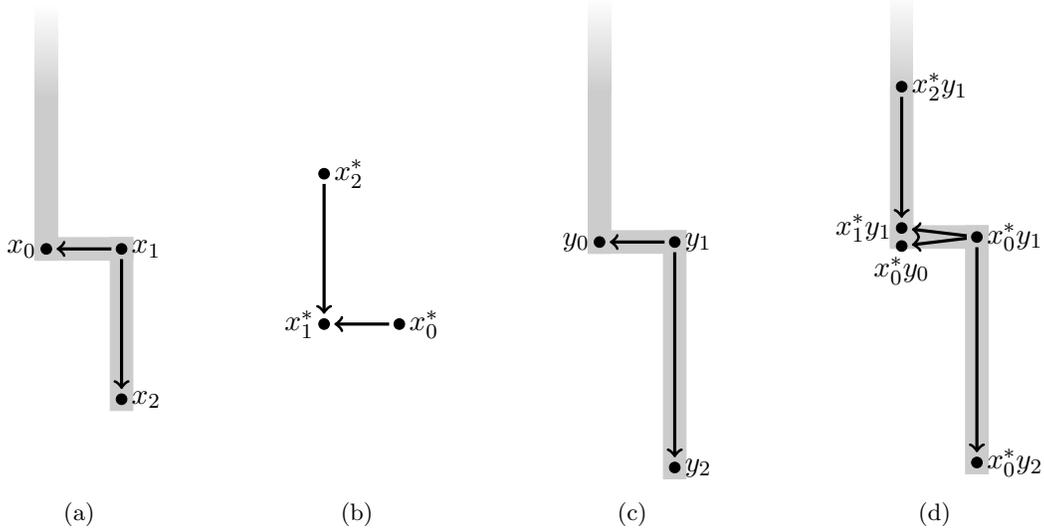

\begin{proof}
We again proceed using induction. We will show that $\varepsilon(J \# -K)=1$ and that
\begin{align*}
a_1(J \# -K) &=a_1(J)\\
a_2(J \# -K) &=a_2(J),
\end{align*}
from which we can conclude that
$$\varepsilon(J \# -nK)=1$$
for all $n \in \N$.

Let $\{x_i\}$ be a basis for $\CFKi(K)$ found using Lemma \ref{lem:basis}. Similarly, let $\{y_i\}$ be such a basis for $\CFKi(J)$. We consider the knot $J \# -K$ and its knot Floer complex. For ease of notation, let $\tau=\tau(J \# -K)$.

Let
\begin{align*}
A&=C\{\textup{min} (i, j-\tau)=0, i \leq a_1(J) \}\\
B&=C\big\{ \{\textup{min} (i, j-\tau)=0, i \leq a_1(J)\} \cup \{ i=a_1(J), \tau-a_2(J) \leq j < \tau \} \big\},
\end{align*}
We claim that the element $x^*_0 y_0 + x^*_1y_1$ generates $H_*(C\{i=0\})$, is zero in $H_*(A)$, and is non-zero in $H_*(B)$. Indeed, there is a direct summand of $B$ with the following generators in the following $(i, j)$-positions:
\begin{align*}
& x^*_0y_0, \ x^*_1y_1 && \big( 0, \tau(J \# -K) \big) \\
& x^*_0y_1 && \big( a_1(J), \tau(J \# -K) \big) \\
& x^*_0y_2 && \big( a_1(J), \tau(J \# -K)-a_2(J) \big)\\
& x^*_2y_1 && \big( 0, \tau(J \# -K) +a_2(K) \big),
\end{align*}
and the following differentials:
\begin{align*}
\partial (x^*_0y_1) &= x^*_0y_0 + x^*_1y_1+x^*_0y_2 \\
\partial (x^*_2y_1)&= x^*_1y_1.
\end{align*}
See Figure \ref{subfig:summand}. From this observation, the claim readily follows; that is,
\begin{align*}
\varepsilon(J \# -K)&=1 \\
a_1(J \# -K) &= a_1(J) \\
a_2(J \# -K) &= a_2(J),
\end{align*}
as desired.
\end{proof}

Recall that an \emph{$L$-space} is a rational homology sphere $Y$ for which
$$\textup{rk } \widehat{HF}(Y) = | H_1(Y, \Z)|.$$
We call a knot $K \subset S^3$ an \emph{$L$-space knot} if there exists $n\in \N$ such that $n$ surgery on $K$ yields an $L$-space. In \cite[Theorem 1.2]{OSlens}, Ozsv\'ath and Szab\'o prove that if $K$ is an $L$-space knot, then the complex $\CFKi(K)$ has a particularly simple form that can be deduced form the Alexander polynomial of $K$, $\Delta_K(t)$. (Note that the results in \cite{OSlens} are stated in terms of $\widehat{HFK}(K)$, but by considering gradings, they are actually sufficient to determine the full $\CFKi(K)$ complex.)

One consequence is that if $K$ is an $L$-space knot, then the Alexander polynomial of $K$ has the form
$$\Delta_K(t) =\sum_{i=0}^k (-1)^i t^{n_i},$$
for some decreasing sequence of non-negative integers $n_0 > n_1 > \ldots > n_k$ with the symmetry condition
$$n_i+n_{k-i}=2g(K),$$
where we have normalized the Alexander polynomial to have a constant term and no negative exponents.
Note that $k$ is always even since there are always an odd number of terms in the Alexander polynomial.

\begin{lemma}
\label{lem:Lspace}
Let $K$ be an $L$-space knot with Alexander polynomial
$$\Delta_K(t) =\sum_{i=0}^k (-1)^i t^{n_i},$$
for some decreasing sequence of integers $n_0 > n_1 > \ldots > n_k$.
Then
\begin{align*}
a_1(K) &=n_0 - n_1 \\
a_2(K) &=n_1-n_2.
\end{align*}
\end{lemma}

\begin{proof}
Theorem 1.2 of  \cite{OSlens} tells us that for $K$ an $L$-space knot, $\widehat{HFK}(K)$ is completely determined by $\Delta_K(t)$. Moreover, up to filtered chain homotopy equivalence, $\CFKi(K)$ is generated as a $\F[U, U^{-1}]$-module by $\widehat{HFK}(K)$, where $\widehat{HFK}$ is the homology of the associated graded object of $\widehat{CFK}(K) \simeq C\{i=0\}$. By considering the gradings on the complex $\CFKi(K)$, and the fact that the differential decreases the Maslov grading by one, the lemma follows.
\end{proof}

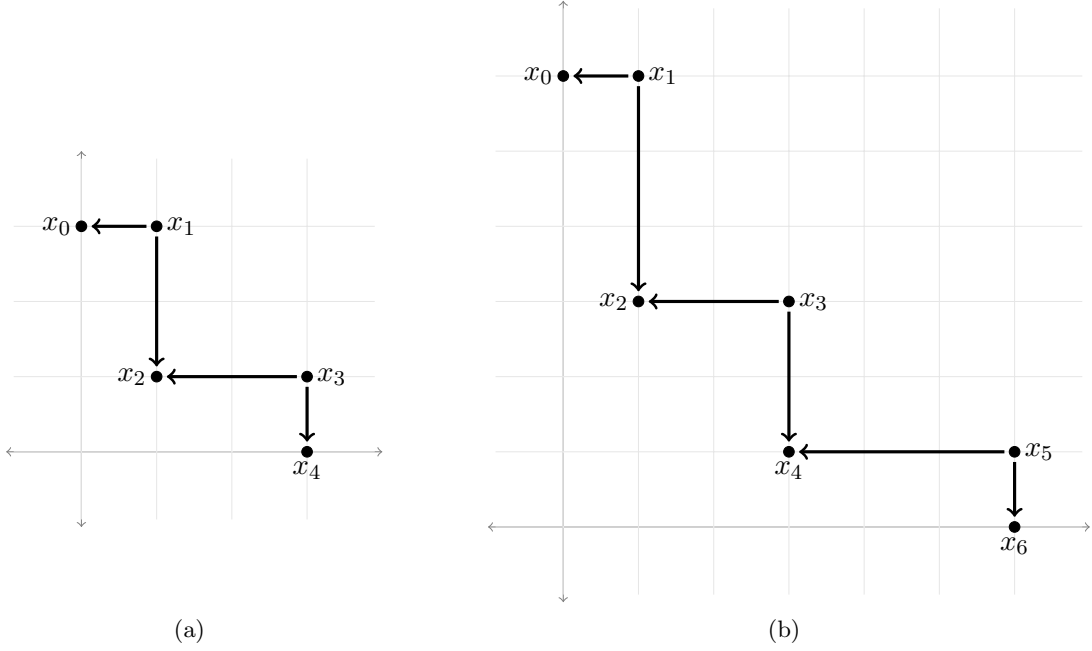
\begin{figure}[htb!]
\vspace{5pt}
\labellist
\small \hair 2pt
\endlabellist
\centering
\subfigure[]{
\begin{tikzpicture}
	\begin{scope}[thin, gray]
		\draw [<->] (-1, 0) -- (4, 0);
		\draw [<->] (0, -1) -- (0, 4);
	\end{scope}
	\draw[step=1, black!10!white, very thin] (-0.9, -0.9) grid (3.9, 3.9);
	\filldraw (0, 3) circle (2pt) node[] (x_0){};
	\filldraw (1, 3) circle (2pt) node[] (x_1){};
	\filldraw (1, 1) circle (2pt) node[] (x_2){};
	\filldraw (3, 1) circle (2pt) node[] (x_3){};
	\filldraw (3, 0) circle (2pt) node[] (x_4){};
	\draw [very thick, <-] (x_0) -- (x_1);
	\draw [very thick, <-] (x_2) -- (x_1);
	\draw [very thick, <-] (x_2) -- (x_3);
	\draw [very thick, <-] (x_4) -- (x_3);
	\node [left] at (x_0) {$x_0$};
	\node [right] at (x_1) {$x_1$};
	\node [left] at (x_2) {$x_2$};
	\node [right] at (x_3) {$x_3$};
	\node [below] at (x_4) {$x_4$};
	\draw (0, -2);
\end{tikzpicture}
}
\hspace{25pt}
\subfigure[]{
\begin{tikzpicture}
	\begin{scope}[thin, gray]
		\draw [<->] (-1, 0) -- (7, 0);
		\draw [<->] (0, -1) -- (0, 7);
	\end{scope}
	\draw[step=1, black!10!white, very thin] (-0.9, -0.9) grid (6.9, 6.9);
	\filldraw (0, 6) circle (2pt) node[] (x_0){};
	\filldraw (1, 6) circle (2pt) node[] (x_1){};
	\filldraw (1, 3) circle (2pt) node[] (x_2){};
	\filldraw (3, 3) circle (2pt) node[] (x_3){};
	\filldraw (3, 1) circle (2pt) node[] (x_4){};
	\filldraw (6, 1) circle (2pt) node[] (x_5){};
	\filldraw (6, 0) circle (2pt) node[] (x_6){};
	\draw [very thick, <-] (x_0) -- (x_1);
	\draw [very thick, <-] (x_2) -- (x_1);
	\draw [very thick, <-] (x_2) -- (x_3);
	\draw [very thick, <-] (x_4) -- (x_3);
	\draw [very thick, <-] (x_4) -- (x_5);
	\draw [very thick, <-] (x_6) -- (x_5);
	\node [left] at (x_0) {$x_0$};
	\node [right] at (x_1) {$x_1$};
	\node [left] at (x_2) {$x_2$};
	\node [right] at (x_3) {$x_3$};
	\node [below] at (x_4) {$x_4$};
	\node [right] at (x_5) {$x_5$};
	\node [below] at (x_6) {$x_6$};
\end{tikzpicture}
}
\caption{Left, the basis from Remark \ref{rem:Lspace} for $\CFKi$ of the torus knot $T_{3, 4}$ with Alexander polynomial $\Delta_{T_{3,4}}(t)=t^6-t^5+t^3-t+1$. Right, the basis for $\CFKi$ of the torus knot $T_{4, 5}$ with Alexander polynomial $\Delta_{T_{4,5}}(t)=t^{12}-t^{11}+t^8-t^6+t^4-t+1$. The lengths of the differentials are given by the differences of the exponents of the Alexander polynomial.}
\label{fig:Lspaceknot}
\end{figure}

\begin{remark}
\label{rem:Lspace}
\emph{
More generally, it can be deduced from \cite[Theorem 1.2]{OSlens} that there is a basis $\{ x_0, \ldots x_k\}$ for $\CFKi(K)$ such that
\begin{align*}
&\partial x_{i} =x_{i-1}+x_{i+1} &&\textup{ for $i$ odd}\\
&\partial x_i =0  &&\textup{ otherwise},
\end{align*}
where the arrow from $x_{i}$ to $x_{i-1}$ is horizontal of length $n_{i}-n_{i-1}$, and the arrow from $x_{i}$ to $x_{i+1}$ is vertical of length $n_{i+1}-n_{i}$. The complex looks like a ``staircase'', where the differences of the $n_i$ give the heights and widths of the steps. See Figure \ref{fig:Lspaceknot}.
}
\end{remark}

Recall that positive torus knots are $L$-space knots since $(pq \pm 1)$-surgery on the torus knot $T_{p,q}$, $p, q>1$, results in a lens space.

\begin{lemma}
\label{lem:Tpp+1}
For $p \geq 3$, the Alexander polynomial of the torus knot $T_{p, p+1}$ is
$$\Delta_{T_{p, p+1}}(t)=\sum_{i=0}^k (-1)^i t^{n_i},$$
for a decreasing sequence of integers $n_0 >n_1> \ldots > n_k$ with
\begin{align*}
n_0 &= p^2-p \\
n_1 &= p^2-p-1 \\
n_2 &= p^2-2p\\
n_3 &= p^2-2p-2.
\end{align*}
In particular,
\begin{align*}
a_1(T_{p, p+1}) &=1 \\
a_2(T_{p, p+1}) &=p-1.
\end{align*}
\end{lemma}

\begin{proof}
Recall that
$$\Delta_{T_{p,q}}(t)=\frac{(t^{pq}-1)(t-1)}{(t^p-1)(t^q-1)}.$$
Following the proof of Proposition 6.1 in \cite{HeddenLivRub}, we see that
$$\frac{(t^{p(p+1)}-1)(t-1)}{(t^p-1)(t^{p+1}-1)} = \sum_{i=0}^{p-1} t^{pi} - t \sum_{i=0}^{p-2} t^{(p+1)i}.$$
Indeed, multiplying both sides by $(t^p-1)(t^{p+1}-1)$, we obtain two telescoping sums on the right-hand side:
\begin{align*}
 (t^p-1)(t^{p+1}-1) \Big( \sum_{i=0}^{p-1} t^{pi} - t \sum_{i=0}^{p-2} t^{(p+1)i} \Big)&= (t^{p+1}-1)(t^{p(p-1)+p}-1) - t(t^p-1)(t^{(p+1)(p-2)+p+1}-1) \\
 &= t^{p^2+p+1} -t^{p^2+p}-t+1 \\
 &= (t^{p(p+1)}-1)(t-1)
\end{align*}
as desired.

The last statement now follows from Lemma \ref{lem:Lspace}.
\end{proof}

\begin{remark}
\emph{For the torus knot $T_{2,3}$, i.e., the case where $p=2$, we can check by hand that
\begin{align*}
a_1(T_{p, p+1}) &=1 \\
a_2(T_{p, p+1}) &=p-1,
\end{align*}
since $\Delta_{T_{2, 3}} (t)=t^2-t+1$.
}
\end{remark}

\begin{remark}
\emph{More generally, for the torus knot $T_{p, p+1}$, the horizontal arrows increase in length by one at each ``step'', from $1$ to $p-1$, and the vertical arrows decrease in length by one at each ``step'', from $p-1$ to $1$. See Figure \ref{fig:Lspaceknot}.}
\end{remark}

\begin{lemma}
\label{lem:iterated}
The iterated torus knot $T_{2,3; p, p+1}$, $p\geq 2$, is an $L$-space knot with Alexander polynomial 
$$\Delta_{T_{2,3; p, p+1}}(t)=\sum_{i=0}^k (-1)^i t^{n_i},$$
for a decreasing sequence of integers $n_0 >n_1> \ldots > n_k$ with
\begin{align*}
n_0 &=p^2+p \\
n_1 &= p^2+p-1 \\
n_2 &= p^2-1.
\end{align*}
In particular,
\begin{align*}
a_1(T_{2,3; p, p+1}) &=1 \\
a_2(T_{2,3; p, p+1}) &=p.
\end{align*}
\end{lemma}

\begin{proof}
The fact that $T_{2,3; \ p, p+1}$ is an $L$-space knot follows from \cite[Theorem 1.10]{HeddencablingII} (cf. \cite{HomLspace}), where Hedden gives sufficient conditions for the cable of an $L$-space knot to again be an $L$-space knot.

The form of the Alexander polynomial follows from the formula for the Alexander polynomial of the cable of knot, i.e.,
$$\Delta_{T_{2,3; p, p+1}}(t)=\Delta_{T_{2,3}}(t^p) \cdot \Delta_{T_{p, p+1}}(t),$$
and Lemma \ref{lem:Tpp+1}. More precisely, for $p \geq 3$,
\begin{align*}
\Delta_{T_{2,3;p, p+1}} &= (t^{2p}-t^p+1)(t^{p^2-p}-t^{p^2-p-1}+t^{p^2-2p}-t^{p^2-2p-2}+ \textup{lower order terms}) \\
&= t^{p^2+p} -t^{p^2+p-1} +t^{p^2-1} + \textup{lower order terms}.
\end{align*}
The case $p=2$ follows easily from the fact that
\[ \Delta_{T_{2,3;2,3}} (t) = t^6 -t^5 + t^3 - t +1.\]
\end{proof}

\begin{lemma}
\label{lem:iterated2}
For $p \geq 2$, $m\geq p^2-p-1$, and $m \neq 1$, the iterated torus knot $T_{p, p+1; 2,2m+1}$ is an $L$-space knot with Alexander polynomial 
$$\Delta_{T_{p, p+1; 2,2m+1}}(t)=\sum_{i=0}^k (-1)^i t^{n_i},$$
for a decreasing sequence of integers $n_0 >n_1> \ldots > n_k$ with
\begin{align*}
n_0 &=2p^2-2p+2m \\
n_1 &= 2p^2-2p+2m-1 \\
n_2 &= 2p^2-4p+2m.
\end{align*}
In particular,
\begin{align*}
a_1(T_{2,3; p, p+1}) &=1 \\
a_2(T_{2,3; p, p+1}) &=2p-1.
\end{align*}
\end{lemma}

\begin{proof}
This iterated torus knot is an $L$-space knot by \cite[Theorem 1.10]{HeddencablingII}. The form of the Alexander polynomial follows from the following facts:
\begin{align*}
\Delta_{T_{p, p+1; 2,2m+1}}(t) &=\Delta_{T_{p, p+1}}(t^2) \cdot \Delta_{T_{2, 2m+1}}(t), \\
\Delta_{T_{2, 2m+1}}(t) &= \sum_{i=0}^{2m} (-1)^i t^i,
\end{align*}
and Lemma \ref{lem:Tpp+1}. More precisely,
\begin{align*}
\Delta_{T_{p, p+1; 2, 2m+1}}(t) &= \Big(  \sum_{i=0}^{p-1} t^{2pi} - t^2 \sum_{i=0}^{p-2} t^{(2p+2)i} \Big) \Big( \sum_{i=0}^{m}  t^{2i} -  \sum_{i=0}^{m-1}  t^{2i+1}\Big) \\
&=  \Big( t^{2p^2-2p} -t^{2p^2-2p-2} + \sum_{i=0}^{p-2} t^{2pi} - t^2 \sum_{i=0}^{p-3} t^{(2p+2)i} \Big) \Big( \sum_{i=0}^{m}  t^{2i} -  \sum_{i=0}^{m-1}  t^{2i+1}\Big) \\
&= t^{2p^2-2p+2m} - t^{2p^2-2p-2} -t^{2p^2-2p+2m-1} +t^{2p^2-2p-1}+ \\ 
& \quad \Big( \sum_{i=0}^{p-2} t^{2pi} - t^2 \sum_{i=0}^{p-3} t^{(2p+2)i} \Big) \Big( \sum_{i=0}^{m}  t^{2i} -  \sum_{i=0}^{m-1}  t^{2i+1}\Big) \\
&= t^{2p^2-2p+2m} - t^{2p^2-2p-2} -t^{2p^2-2p+2m-1} +t^{2p^2-2p-1}+  \\
& \quad \big(t^{2p^2-4p+2m} +\textup{lower order terms} \big) \\ 
&=  t^{2p^2-2p+2m} -t^{2p^2-2p+2m-1} + t^{2p^2-4p+2m} +\textup{lower order terms},
\end{align*}
where the last equality follows from the hypothesis that $m>p$.
\end{proof}

Recall that $D$ denotes the (positive, untwisted) Whitehead double of the right-handed trefoil.
\begin{lemma}
\label{lem:D}
As elements of the group $\cF$,
$$\llbracket D\rrbracket =\llbracket T_{2,3}\rrbracket .$$
\end{lemma}

\begin{proof}
In \cite[Theorem 1.2]{HeddenWhitehead}, Hedden determines the $\Z$-filtered chain homotopy type of $\widehat{CFK}$ of the Whitehead double of $K$ in terms of $\widehat{CFK}(K)$. We can use this result to determine $\widehat{CFK}(D)$, from which we will deduce the class $\llbracket D\rrbracket $ using rank and grading considerations.

Using Hedden's result, we see that
$$\widehat{CFK}(D, j) \simeq \left\{
	\begin{array}{ll}
		\F^2_{(0)} \oplus \F^2_{(-1)} & \quad j=1 \\
		\F^3_{(-1)} \oplus \F^4_{(-2)} & \quad j=0 \\
		\F^2_{(-2)} \oplus  \F^2_{(-3)} & \quad j=-1
	\end{array} \right.$$
where the subscript denotes the Maslov, or homological, grading, and $j$ denotes the Alexander grading. Moreover, Hedden proves that every non-trivial differential on this complex lowers the Alexander grading by exactly one, which is sufficient to completely determine the $\Z$-filtered chain homotopy type of $\widehat{CFK}(D)$. Note that $\tau(D)=1$.

Let $x$ be a generator of $\widehat{HF}(S^3) \cong H_*(C\{i=0\})$. Note that $x$ necessarily is positioned at $(0, 1)$ in the $(i,j)$-plane. Then $[x]$ must be zero in $H_*(C\{j=1\})$ since the homology of $C\{ j=1 \}$ is supported in $i$-coordinate $2$. By considering the support of $\widehat{CFK}(D)$, we see that $x$ is in the kernel of $\partial^\horz$, so in order to vanish in $H_*(C\{j=1\})$, it must be in the image of $\partial^\horz$, i.e., there exists a class, say $y$, positioned at $(1,1)$, such that
$$\partial^\horz y = x.$$
The class $[y]$ is equal to zero in $H_*(C\{i=1\})$ since the homology of $C\{i=1\}$ is supported in $j$-coordinate $2$. But $y$ cannot be in the image of the differential on $C\{i=1\}$, since $\partial^2=0$, where $\partial$ is the differential on $\CFKi$, and $\partial^\horz y \neq 0$. Hence, the boundary of $y$ in $C\{i=1\}$ must be non-zero; denote this boundary by $z$. Notice that $z$ has $(i, j)$-coordinates $(1, 0)$.

Again, for $\partial^2=0$ reasons, the boundary of $z$ in $C\{j=0\}$ must be zero, and by grading considerations, $z$ is not in the image of the differential on $C\{j=0\}$.

The complex $\CFKi(-T_{2,3})$ is generated over $\F[U, U^{-1}]$ by
$$a, \ b, \ c,$$
with the differential
\begin{align*}
\partial a &=b \\
\partial c &= b,
\end{align*}
where the generators are have the following $(i, j)$-coordinates:
\begin{align*}
a \qquad & (0, 1) \\
b \qquad & (0, 0) \\
c \qquad  & (1, 0).
\end{align*}
Then in the tensor product
$$\CFKi(-T_{2, 3}) \otimes_{\F[U, U^{-1}]} \CFKi(D)$$
the generator
$$az+by+cx$$
is non-trivial in both vertical and horizontal homology. Indeed, it is clearly in the kernel of the vertical differential, and cannot be in the image of the vertical differential, since $cx$ does not appear in the vertical boundary of any element. Similarly, it is in the kernel but not the image of the horizontal differential.

Thus,
$$\varepsilon\big(\CFKi(-T_{2, 3}) \otimes_{\F[U, U^{-1}]} \CFKi(D) \big) =0,$$
as desired.
\end{proof}

We are now ready to prove Proposition \ref{prop:examples}, showing that we have the following relations in $\cF$, where $0<p<q$:
\begin{itemize}
	\item $\llbracket T_{p, p+1}\rrbracket  \ll \llbracket T_{q, q+1}\rrbracket $
	\item \label{it:cableofdouble} $\llbracket D_{p, p+1}\rrbracket  \ll \llbracket D_{q, q+1}\rrbracket $
	\item \label{it:torusdouble} $\llbracket T_{p, p+1}\rrbracket  \ll \llbracket D_{p, p+1}\rrbracket $
	\item $\llbracket T_{p, p+1; 2, 2m+1}\rrbracket  \ll \llbracket T_{q, q+1; 2, 2m+1}\rrbracket $, for $m \geq q^2-q-1$.
\end{itemize}

\begin{proof}[Proof of Proposition \ref{prop:examples}]
The proposition is now an easy consequence of the preceding lemmas. We have from Lemma \ref{lem:Tpp+1} that
\begin{align*}
a_1(T_{p, p+1}) &=1 \\
a_2(T_{p, p+1}) &=p-1.
\end{align*}
Now Lemma \ref{lem:agg} states that if $a_1(J)=a_1(K)$ and $a_2(J) <a_2(K)$, then $\llbracket J\rrbracket  \ll \llbracket K\rrbracket $, implying that
$$\llbracket T_{p, p+1}\rrbracket  \ll \llbracket T_{q, q+1}\rrbracket ,$$
which proves the first assertion in the proposition.

From Lemma \ref{lem:D}, we have that
$$\llbracket D\rrbracket =\llbracket T_{2,3}\rrbracket ,$$
and from Proposition \ref{prop:satellitemap} that
$$\llbracket D_{p, p+1}\rrbracket =\llbracket T_{2,3; p, p+1}\rrbracket .$$
Hence by Lemmas \ref{lem:class} and \ref{lem:iterated},
\begin{align*}
a_1(D_{p, p+1}) &=1 \\
a_2(D_{p, p+1}) &=p,
\end{align*}
so by Lemma \ref{lem:agg},
$$\llbracket D_{p, p+1}\rrbracket  \ll \llbracket D_{q, q+1}\rrbracket  \quad \textup{ and } \quad \llbracket T_{p, p+1}\rrbracket  \ll \llbracket D_{p, p+1}\rrbracket .$$
Finally, by Lemma \ref{lem:iterated2}, we have that 
\begin{align*}
a_1(T_{p, p+1;2, 2m+1}) &=1 \\
a_2(T_{p, p+1;2, 2m+1}) &= 2p-1,
\end{align*}
for $p\geq 2$, $m \geq p^2-p-1$, $m \neq 1$, and so
$$\llbracket T_{p, p+1; 2, 2m+1} \rrbracket \ll \llbracket T_{q, q+1; 2, 2m+1} \rrbracket.$$
This completes the proof of the proposition.
\end{proof}

We conclude this paper by showing that our examples, $\{D_{p, p+1} \# -T_{p, p+1}\}_{p \geq 2}$, of smoothly independent, topologically slice knots are smoothly independent  from the examples of Endo \cite{Endo} and Hedden-Kirk \cite{HeddenKirk}. Recall that Endo's examples are pretzel knots of the form
$$K_t=K(-2t-1, 4t+1, 4t+3), \quad t \geq 1.$$
In particular, they are of genus one.
The examples of Hedden-Kirk are (positive, untwisted) Whitehead doubles of certain torus knots.

\begin{proposition} \label{prop:genusone}
If $K$ is a knot of genus one and $\varepsilon(K)=1$, then either
$$a_1(K) \neq 1 \qquad \textup{ or } \qquad a_1(K)=a_2(K)=1.$$
\end{proposition}

\begin{proof}
Notice that the assumption that $\varepsilon(K)=1$ does not cause any loss of generality, since $\varepsilon(-K)=-\varepsilon(K)$.

Assume that $a_1(K) =1$. We first notice that if $K$ is a knot of genus one and $\varepsilon(K)=1$, then $\tau(K) \neq -1$. This follows from the adjunction inequality for knot Floer homology \cite[Theorem 5.1]{OSknots}, and the basis from Lemma \ref{lem:basis}

Now, suppose $a_1(K)=1$ and $\tau(K)=0$. Using the adjunction inequality \cite[Theorem 5.1]{OSknots}, and a basis found using the first part of Lemma \ref{lem:basis}, we see that the basis element $x_1$ must be in the kernel of the differential on $C\{i=1\}$. Moreover, for $\partial^2=0$ reasons, it cannot be in the image of the differential on $C\{i=1\}$. But $[x_1]$ cannot be zero in $H_*(C\{i=1\}$, because $\tau(K)=0$ implies that $H_*(C\{i=1\})$ is supported in $(i, j)$-coordinate $(1,1)$. 
 
Hence, we may assume that $a_1(K)=1$ and $\tau(K)=1$, in which case the arguments in the proof of Lemma \ref{lem:D} lead us to the desired result.
\end{proof}

\noindent In the proof of Proposition \ref{prop:examples}, we showed that
\begin{align*}
a_1(D_{p, p+1} \# -T_{p, p+1}) &=1 \\
a_2(D_{p, p+1} \# -T_{p, p+1}) &=p,
\end{align*}
Hence, by Proposition \ref{prop:genusone}, along with Lemmas \ref{lem:agg1} and \ref{lem:agg}, it follows that when $p>1$, our examples are independent from those of Endo and Hedden-Kirk.

The following proposition describes the subgroup of $\cF$ generated by Whitehead doubles:
\begin{proposition}
Whitehead doubles are contained in the rank one subgroup of $\cF$ generated by the right-handed trefoil.
\end{proposition}

\begin{proof}
The argument in Lemma \ref{lem:D} can be used to show that for a Whitehead double $WD$ with $\varepsilon(WD)=1$, the class $\llbracket WD\rrbracket =\llbracket T_{2,3}\rrbracket $ in $\cF$. This is sufficient for the result, since $\varepsilon(WD)=-1$ implies that $\varepsilon(-WD)=1$, and $\varepsilon(WD)=0$ implies that $\llbracket WD \rrbracket = 0$.
\end{proof}

\bibliographystyle{amsalpha}

\bibliography{mybib}

\end{document}